\documentclass[11pt]{article}

\usepackage{latexsym}
\usepackage[dvips]{graphicx}
\usepackage{float}
\usepackage{amsmath}
\usepackage{amssymb}
\usepackage{epic}
\usepackage{color}
\usepackage{setspace}
\usepackage{lscape}
\usepackage{arydshln}
\usepackage{fancybox}   
\usepackage{mathtools}
\usepackage{amsthm}

\begin{document}
\bibliographystyle{plain}
\floatplacement{table}{H}
\newtheorem{definition}{Definition}[section]
\newtheorem{lemma}{Lemma}[section]
\newtheorem{theorem}{Theorem}[section]
\newtheorem{corollary}{Corollary}[section]
\newtheorem{proposition}{Proposition}[section]

\newcommand{\sni}{\sum_{i=1}^{n}}
\newcommand{\snj}{\sum_{j=1}^{n}}
\newcommand{\smj}{\sum_{j=1}^{m}}
\newcommand{\sumjm}{\sum_{j=1}^{m}}
\newcommand{\bdis}{\begin{displaymath}}
\newcommand{\edis}{\end{displaymath}}
\newcommand{\beq}{\begin{equation}}
\newcommand{\eeq}{\end{equation}}
\newcommand{\beqn}{\begin{eqnarray}}
\newcommand{\eeqn}{\end{eqnarray}}
\newcommand{\defeq}{\stackrel{\triangle}{=}}
\newcommand{\simleq}{\stackrel{<}{\sim}}
\newcommand{\sep}{\;\;\;\;\;\; ; \;\;\;\;\;\;}
\newcommand{\real}{\mbox{$ I \hskip -4.0pt R $}}
\newcommand{\complex}{\mbox{$ I \hskip -6.8pt C $}}
\newcommand{\integ}{\mbox{$ Z $}}
\newcommand{\realn}{\real ^{n}}
\newcommand{\sqrn}{\sqrt{n}}
\newcommand{\sqrtwo}{\sqrt{2}}
\newcommand{\prf}{{\bf Proof. }}

\newcommand{\onehlf}{\frac{1}{2}}
\newcommand{\thrhlf}{\frac{3}{2}}
\newcommand{\fivhlf}{\frac{5}{2}}
\newcommand{\onethd}{\frac{1}{3}}
\newcommand{\lb}{\left ( }
\newcommand{\lcb}{\left \{ }
\newcommand{\lsb}{\left [ }
\newcommand{\labs}{\left | }
\newcommand{\rb}{\right ) }
\newcommand{\rcb}{\right \} }
\newcommand{\rsb}{\right ] }
\newcommand{\rabs}{\right | }
\newcommand{\lnm}{\left \| }
\newcommand{\rnm}{\right \| }
\newcommand{\lambdab}{\bar{\lambda}}
%
%           I N D I C E S
%
\newcommand{\xj}{x_{j}}
\newcommand{\xjb}{\bar{x}_{j}}
\newcommand{\xro}{x_{\resh}}
\newcommand{\xrob}{\bar{x}_{\resh}}
\newcommand{\xsig}{x_{\sigma}}
\newcommand{\xsigb}{\bar{x}_{\sigma}}
\newcommand{\xnmjb}{\bar{x}_{n-j+1}}
\newcommand{\xnmj}{x_{n-j+1}}
\newcommand{\aroj}{a_{\resh j}}
\newcommand{\arojb}{\bar{a}_{\resh j}}
\newcommand{\aroro}{a_{\resh \resh}}
\newcommand{\amuro}{a_{\mu \resh}}
\newcommand{\amumu}{a_{\mu \mu}}
\newcommand{\aii}{a_{ii}}
\newcommand{\aik}{a_{ik}}
\newcommand{\akj}{a_{kj}}
\newcommand{\atwoii}{a^{(2)}_{ii}}
\newcommand{\atwoij}{a^{(2)}_{ij}}
\newcommand{\ajj}{a_{jj}}
\newcommand{\aiib}{\bar{a}_{ii}}
\newcommand{\ajjb}{\bar{a}_{jj}}
\newcommand{\bii}{a_{jj}}
\newcommand{\biib}{\bar{a}_{jj}}
\newcommand{\aij}{a_{i,n-i+1}}
\newcommand{\akl}{a_{j,n-j+1}}
\newcommand{\aijb}{\bar{a}_{i,n-i+1}}
\newcommand{\aklb}{\bar{a}_{j,n-j+1}}
\newcommand{\bij}{a_{n-j+1,j}}
\newcommand{\arorob}{\bar{a}_{\resh \resh}}
\newcommand{\arosig}{a_{\resh \sigma}}
\newcommand{\arosigb}{\bar{a}_{\resh \sigma}}
\newcommand{\sumjrosig}{\sum_{\stackrel{j=1}{j\neq\resh,\sigma}}^{n}}
\newcommand{\summuro}{\sum_{\stackrel{j=1}{j\neq\mu,\resh}}^{n}}
\newcommand{\sumjnoti}{\sum_{\stackrel{j=1}{j\neq i}}^{n}}
\newcommand{\sumlnoti}{\sum_{\stackrel{\ell=1}{\ell \neq i}}^{n}}
\newcommand{\sumknoti}{\sum_{\stackrel{k=1}{k\neq i}}^{n}}
\newcommand{\sumknotij}{\sum_{\stackrel{k=1}{k\neq i,j}}^{n}}
\newcommand{\sumk}{\sum_{k=1}^{n}}
\newcommand{\snl}{\sum_{\ell=1}^{n}}
\newcommand{\sumji}{\sum_{\stackrel{j=1}{j\neq i, n-i+1}}^{n}}
\newcommand{\sumki}{\sum_{\stackrel{k=1}{k\neq i, n-i+1}}^{n}}
\newcommand{\sumkj}{\sum_{\stackrel{k=1}{k\neq j, n-j+1}}^{n}}
\newcommand{\sumjro}{\sum_{\stackrel{j=1}{j\neq\resh}}^{n}}
\newcommand{\rrosig}{R''_{\resh \sigma}}
\newcommand{\rro}{R'_{\resh}}
\newcommand{\gamror}{\Gamma_{\resh}^{R}(A)}
\newcommand{\gamir}{\Gamma_{i}^{R}(A)}
\newcommand{\gamctrr}{\Gamma_{\frac{n+1}{2}}^{R}(A)}
\newcommand{\gamctrc}{\Gamma_{\frac{n+1}{2}}^{C}(A)}
\newcommand{\gamroc}{\Gamma_{\resh}^{C}(A)}
\newcommand{\gamjc}{\Gamma_{j}^{C}(A)}
\newcommand{\lamror}{\Lambda_{\resh}^{R}(A)}
\newcommand{\lamir}{\Lambda_{i}^{R}(A)}
\newcommand{\lamirepsilon}{\Lambda_{i}^{R}(A_{\epsilon})}
\newcommand{\lamnir}{\Lambda_{n-i+1}^{R}(A)}
\newcommand{\lamjr}{\Lambda_{j}^{R}(A)}
\newcommand{\varphiij}{\Phi_{ij}^{R}(A)}
\newcommand{\delir}{\Delta_{i}^{R}(A)}
\newcommand{\vir}{V_{i}^{R}(A)}
\newcommand{\pamir}{\Pi_{i}^{R}(A)}
\newcommand{\xir}{\Xi_{i}^{R}(A)}
\newcommand{\lamjc}{\Lambda_{j}^{C}(A)}
\newcommand{\vjc}{V_{j}^{C}(A)}
\newcommand{\pamjc}{\Pi_{j}^{C}(A)}
\newcommand{\xjc}{\Xi_{j}^{C}(A)}
\newcommand{\lamroc}{\Lambda_{\resh}^{C}(A)}
\newcommand{\lamsigr}{\Lambda_{\sigma}^{R}(A)}
\newcommand{\lamsigc}{\Lambda_{\sigma}^{C}(A)}
\newcommand{\psii}{\Psi_{i}^{R}(A)}
\newcommand{\psiq}{\Psi_{q}}
\newcommand{\psiiepsilon}{\Psi_{i}^{R}(A_{\epsilon})}
\newcommand{\psiqepsilon}{\Psi_{q}(A_{\epsilon})}
\newcommand{\psiqc}{\Psi_{q}^{c}}
\newcommand{\psiqcepsilon}{\Psi_{q}^{c}(A_{\epsilon})}

\newcommand{\xmu}{x_{\mu}}
\newcommand{\xmub}{\bar{x}_{\mu}}
\newcommand{\xnu}{x_{\nu}}
\newcommand{\xnub}{\bar{x}_{\nu}}
\newcommand{\amuj}{a_{\mu j}}
\newcommand{\amujb}{\bar{a}_{\mu j}}
\newcommand{\amumub}{\bar{a}_{\mu \mu}}
\newcommand{\amunu}{a_{\mu \nu}}
\newcommand{\amunub}{\bar{a}_{\mu \nu}}
\newcommand{\sumjmunu}{\sum_{\stackrel{j=1}{j\neq\mu,\nu}}}
\newcommand{\rmunu}{R''_{\mu \nu}}
\newcommand{\rmu}{R'_{\mu}}

\newcommand{\Azero}{A_{0}}
\newcommand{\Aone}{A_{1}}
\newcommand{\Atwo}{A_{2}}
\newcommand{\Ath}{A_{3}}
\newcommand{\Afr}{A_{4}}
\newcommand{\Afv}{A_{5}}
\newcommand{\Asx}{A_{6}}
\newcommand{\Anmo}{A_{n-1}}
\newcommand{\Anmt}{A_{n-2}}
\newcommand{\An}{A_{n}}
\newcommand{\Aj}{A_{j}}

\newcommand{\azero}{a_{0}}
\newcommand{\aone}{a_{1}}
\newcommand{\atwo}{a_{2}}
\newcommand{\ath}{a_{3}}
\newcommand{\afr}{a_{4}}
\newcommand{\afv}{a_{5}}
\newcommand{\asx}{a_{6}}
\newcommand{\anmo}{a_{n-1}}
\newcommand{\anmt}{a_{n-2}}
\newcommand{\an}{a_{n}}
\newcommand{\aj}{a_{j}}

\newcommand{\qzero}{q_{0}}
\newcommand{\qth}{q_{3}}
\newcommand{\qfr}{q_{4}}
\newcommand{\qfv}{q_{5}}
\newcommand{\qsx}{q_{6}}
\newcommand{\qnmo}{q_{n-1}}
\newcommand{\qnmt}{q_{n-2}}
\newcommand{\qn}{q_{n}}
\newcommand{\qj}{q_{j}}

\newcommand{\Bzero}{B_{0}}
\newcommand{\Bone}{B_{1}}
\newcommand{\Btwo}{B_{2}}
\newcommand{\Bth}{B_{3}}
\newcommand{\Bfr}{B_{4}}
\newcommand{\Bfv}{B_{5}}
\newcommand{\Bsx}{B_{6}}
\newcommand{\Bnmo}{B_{n-1}}
\newcommand{\Bndtwomo}{B_{n/2-1}}
\newcommand{\Bnmt}{B_{n-2}}
\newcommand{\Bn}{B_{n}}
\newcommand{\Bj}{B_{j}}

\newcommand{\Anmk}{A_{n-k}}
\newcommand{\Bnmk}{A_{n-k}}
\newcommand{\anmk}{a_{n-k}}
\newcommand{\zj}{z^{j}}          
\newcommand{\zk}{z^{k}}          

\newcommand{\cii}{c_{ii}}
\newcommand{\cik}{c_{ik}}
\newcommand{\ckj}{c_{kj}}
\newcommand{\ctwoii}{c^{(2)}_{ii}}
\newcommand{\ctwoij}{c^{(2)}_{ij}}
\newcommand{\cjj}{c_{jj}}

\newcommand{\bik}{b_{ik}}
\newcommand{\bkj}{b_{kj}}
\newcommand{\btwoii}{b^{(2)}_{ii}}
\newcommand{\btwoij}{b^{(2)}_{ij}}
\newcommand{\bjj}{b_{jj}}

\newcommand{\abii}{(AB)_{ii}}
\newcommand{\abil}{(AB)_{i\ell}}

\newcommand{\bkl}{b_{k\ell}}
\newcommand{\btwoil}{b^{(2)}_{i\ell}}
\newcommand{\bll}{b_{\ell \ell}}

\newcommand{\matrixspace}{\;\;}
\newcommand{\ellone}{\ell_{1}}
\newcommand{\elltwo}{\ell_{2}}

\newcommand{\varphik}{\varphi_{k}}
\newcommand{\chik}{\chi_{k}}
\newcommand{\Phik}{\Phi_{k}}
\newcommand{\psik}{\psi_{k}}
\newcommand{\dr}{\beta^{n-k}}
\newcommand{\dn}{\beta^{-k}}
\newcommand{\betatwotilde}{\tilde{\rtwonk}}
\newcommand{\betathtilde}{\tilde{\ronenk}}
\newcommand{\betamink}{\beta^{-k}}

\newcommand{\rvarphik}{r(\varphi_{k})}    
\newcommand{\rrvarphik}{(r(\varphi_{k}))}    
\newcommand{\rchik}{r(\chi_{k})}    
\newcommand{\rrchik}{(r(\chi_{k}))}    
\newcommand{\svarphik}{s(\varphi_{k})}    
\newcommand{\ssvarphik}{(s(\varphi_{k}))}    
\newcommand{\schik}{s(\chi_{k})}    
\newcommand{\sschik}{(s(\chi_{k}))}    

\newcommand{\rpsik}{r(\psi_{k})}    
\newcommand{\rrpsik}{(r(\psi_{k}))} 
\newcommand{\spsik}{s(\psi_{k})}    
\newcommand{\sspsik}{(s(\psi_{k}))} 
\newcommand{\rpsinmk}{r(\psi_{n-k})}    
\newcommand{\rrpsinmk}{(r(\psi_{n-k}))} 
\newcommand{\spsinmk}{s(\psi_{n-k})}    
\newcommand{\sspsinmk}{(s(\psi_{n-k}))} 

\newcommand{\rvarphinmk}{r(\varphi_{n-k})}    
\newcommand{\rrvarphinmk}{(r(\varphi_{n-k}))}    
\newcommand{\rchinmk}{r(\chi_{n-k})}    
\newcommand{\rrchinmk}{(r(\chi_{n-k}))}    
\newcommand{\svarphinmk}{s(\varphi_{n-k})}    
\newcommand{\ssvarphinmk}{(s(\varphi_{n-k}))}    
\newcommand{\schinmk}{s(\chi_{n-k})}    
\newcommand{\sschinmk}{(s(\chi_{n-k}))}    
\newcommand{\stilde}{\tilde{s}}

\newcommand{\resh}{\rho}
\newcommand{\snk}{s}
\newcommand{\ronenk}{r_{1}}
\newcommand{\rtwonk}{r_{2}}
\newcommand{\obar}{\widebar{O}}
\newcommand{\Lbar}{\widebar{L}}
\newcommand{\dltaone}{\delta_{1}}
\newcommand{\dltatwo}{\delta_{2}}
\newcommand{\dltatld}{\tilde{\delta}}
\newcommand{\tauone}{\tau_{1}}
\newcommand{\tautwo}{\tau_{2}}
\newcommand{\xb}{\bar{x}}
\newcommand{\qone}{q_{1}}
\newcommand{\qtwo}{q_{2}}
\newcommand{\ub}{\bar{u}}
\newcommand{\cmm}{\complex^{m \times m}}
\newcommand{\opdsk}{\mathcal{O}}
\newcommand{\cldsk}{\widebar{\mathcal{O}}}

\newcommand{\sqrtaone}{\sqrt{\a1}}
\newcommand{\sqrtatwo}{\sqrt{\atwo}}

\newcommand{\rhoone}{\rho_{1}}
\newcommand{\rhotwo}{\rho_{2}}
\newcommand{\xone}{x_{1}}
\newcommand{\xtwo}{x_{2}}
\newcommand{\xthr}{x_{3}}
\newcommand{\xfor}{x_{4}}
\newcommand{\xfiv}{x_{5}}

\newcommand{\lnorm}{\left \|}
\newcommand{\rnorm}{\right \|}
\newcommand{\lnrm}{\biggl | \biggl |}
\newcommand{\rnrm}{\biggr |\biggr |}

\newcommand{\recipp}{p^{\protect \#}}
\newcommand{\recipP}{P^{\protect \#}}

\newcommand{\varphione}{\varphi_{1}}
\newcommand{\varphitwo}{\varphi_{2}}

\newcommand{\psione}{\psi_{1}}
\newcommand{\psitwo}{\psi_{2}}

\newcommand{\absatwo}{|\atwo|}
\newcommand{\absqrtatwo}{\sqrt{|\atwo|}}
\newcommand{\absrhoone}{|\rhoone|}
\newcommand{\absrhotwo}{|\rhotwo|}
\newcommand{\sqrgamma}{\sqrt{\gamma}}
\newcommand{\aminb}{a-b}                 
\newcommand{\sqrtc}{\sqrt{c}}            
\newcommand{\sqrtabsc}{\sqrt{|c|}}            
\newcommand{\absalf}{|\alpha|}            
\newcommand{\rhoi}{\rho_{i}}                  
\newcommand{\sigmai}{\sigma_{i}}                  
\newcommand{\rhoj}{\rho_{j}}                  
\newcommand{\sigmaj}{\sigma_{j}}                  
\newcommand{\xij}{x_{ij}}                  
\newcommand{\yij}{y_{ij}}                  
\newcommand{\varphii}{\varphi_{i}}
\newcommand{\varphij}{\varphi_{j}}
\newcommand{\xstar}{x^{\ast}}
\newcommand{\uuone}{\mathcal{U}_{1}}
\newcommand{\uutwo}{\mathcal{U}_{2}}
\newcommand{\xkova}{\hat{x}}
\newcommand{\xkovalam}{\xkova_{max}(\ell-1)}

\newcommand{\dnz}{(d_{n}(z))}
\newcommand{\comm}[1]{}

\newcommand{\fzero}{f_{0}}
\newcommand{\fone}{f_{1}}
\newcommand{\ftwo}{f_{2}}
\newcommand{\fzeroz}{f_{0}(z)}
\newcommand{\fonez}{f_{1}(z)}
\newcommand{\ftwoz}{f_{2}(z)}
\newcommand{\qzeroz}{q_{0}(z)}
\newcommand{\qonez}{q_{1}(z)}
\newcommand{\qtwoz}{q_{2}(z)}

\newcommand{\otms}{\otimes}

\begin{center}
\large
{\bf EIGENVALUE BOUNDS FOR MATRIX POLYNOMIALS IN GENERALIZED BASES}
\vskip 0.5cm
\normalsize
A. Melman \\
Department of Applied Mathematics \\
School of Engineering, Santa Clara University  \\
Santa Clara, CA 95053  \\
e-mail : amelman@scu.edu \\
\vskip 0.5cm
\end{center}

\begin{abstract}
We derive inclusion regions for the eigenvalues of matrix polynomials expressed in a general polynomial basis, which can lead to significantly
better results than traditional bounds. We present several applications to engineering problems.
\vskip 0.15cm
{\bf AMS(MOS) subject classification :} 12D10, 15A18, 30C15
\end{abstract}

%%%%%%%%%%%%%%%%%%%%%%%%%%%%%%%%%%%%%%%%%%%%%%%%%%%%%%%%%%%%%%%%%%%%%%%%%%%%%%%%%%%%%%%%%%%%%%%%%%%%%%%%%%%%%%%%%%%%%%%%%%%%%%%%%%%%%%%%%%%%%%%%%%%%%%%%%%%%%%%%%%
%%%%%%%%%%%%%%%%%%%%%%%%%%%%%%%%%%%%%%%%%%%%%%%%%%%%%%%%%%%%%%%%%%%%%%%%%%%%%%%%%%%%%%%%%%%%%%%%%%%%%%%%%%%%%%%%%%%%%%%%%%%%%%%%%%%%%%%%%%%%%%%%%%%%%%%%%%%%%%%%%%
%%%%%%%%%%%%%%%%%%%%%%%%%%%%%%%%%%%%%%%%%%%%%%%%%%%%%%%%%%%%%%%%%%%%%%%%%%%%%%%%%%%%%%%%%%%%%%%%%%%%%%%%%%%%%%%%%%%%%%%%%%%%%%%%%%%%%%%%%%%%%%%%%%%%%%%%%%%%%%%%%%
%
%
%              SECTION 1  -  INTRODUCTION                                              
%
%
%%%%%%%%%%%%%%%%%%%%%%%%%%%%%%%%%%%%%%%%%%%%%%%%%%%%%%%%%%%%%%%%%%%%%%%%%%%%%%%%%%%%%%%%%%%%%%%%%%%%%%%%%%%%%%%%%%%%%%%%%%%%%%%%%%%%%%%%%%%%%%%%%%%%%%%%%%%%%%%%%%
%%%%%%%%%%%%%%%%%%%%%%%%%%%%%%%%%%%%%%%%%%%%%%%%%%%%%%%%%%%%%%%%%%%%%%%%%%%%%%%%%%%%%%%%%%%%%%%%%%%%%%%%%%%%%%%%%%%%%%%%%%%%%%%%%%%%%%%%%%%%%%%%%%%%%%%%%%%%%%%%%%
%%%%%%%%%%%%%%%%%%%%%%%%%%%%%%%%%%%%%%%%%%%%%%%%%%%%%%%%%%%%%%%%%%%%%%%%%%%%%%%%%%%%%%%%%%%%%%%%%%%%%%%%%%%%%%%%%%%%%%%%%%%%%%%%%%%%%%%%%%%%%%%%%%%%%%%%%%%%%%%%%%

\section{Introduction}           
\label{introduction} 

A polynomial eigenvalue problem consists in computing a nonzero complex eigenvector $v$ and a complex eigenvalue $z$ such that $P(z)v=0$, 
where $P$ is a matrix polynomial 
of the form
\bdis
\An  z^{n} + \Anmo z^{n-1} + \dots \Aone z + \Azero  \; ,
\edis
and $\Aj$ ($j=1,\dots,n$) are complex $m \times m$ matrices. If $A_{n}$ is singular then there are infinite eigenvalues, and if $A_{0}$ is singular then zero 
is an eigenvalue. There are $nm$ eigenvalues, including possibly infinite ones. The finite eigenvalues are the solutions 
of $\text{det} P(z)=0$. We refer to~\cite{BHMST} and~\cite{TisseurMeerbergen} for an overview of engineering applications. 

It is, in general, a computationally intensive task to solve these problems, although bounds on the eigenvalues are relatively easy to compute. 
Such bounds are useful, e.g., in eigenvalue computation by iterative methods (\cite{SP}) and when computing pseudospectra 
(\cite{HT_pseudospectra},\cite{TH_pseudospectra}). Most localization results for polynomial eigenvalues found in the literature apply to
matrix polynomials that are expressed in the regular polynomial power basis $\{1,z,z^{2}, \dots\}$. However, using a different basis can lead to significantly
better results for particular classes of problems. It is what we propose to do here.

Bounds for matrix polynomials are often based on bounds for scalar polynomials, many examples of which can be found in~\cite{HighamTisseur}, and our 
approach will be similar. Specifically, we were inspired by Theorem~8.4.6 in~\cite{RS}, where a zero inclusion region is derived for
scalar polynomials, expressed in a weakly interlacing basis, namely, a basis consisting of polynomials with real zeros that weakly interlace. 
We found such bases to be too restrictive, and we will derive a matrix version of a generalization of this theorem to more general bases. 
For one of those bases, the Newton basis with complex nodes, our result, applied to scalar polynomials, is Theorem~8.6.3 of~\cite{RS}.
We therefore unify and generalize to matrix polynomials inclusion regions not only for weakly interlacing and Newton bases, but for more general ones as well.

To make our exposition reasonably self-contained, we now state a few theorems and definitions that we will use later.
The first is an extension to matrix-valued analytical functions of Rouch\'{e}'s theorem from~\cite{BiniNoferiniSharify} and~\cite{Melman_MatPol}.
Throughout, all matrix norms will assumed to be subordinate, i.e., induced by a vector norm.
%%%%%%%%%%%%%%%%%%%%%%%%%%%%%%%%%%%%%%%%%%%%%%%%%%%%%%%%%%%%%%%%%%%%%%%%%%%%%%%%%%%%%%%%%%%%%%%%%%%%%%%%%%%%%%%%%%%%%%%%%%%%%%%%%%%%%%%%%%%%%%%%%%%%%%%%%%%%%%%%%%
%%%%%%%%%%%%%%%%%%%%%%%%%%%%%%%%%%%%%%%%%%%%%%%%%%%%%%%%%%%%%%%%%%%%%%%%%%%%%%%%%%%%%%%%%%%%%%%%%%%%%%%%%%%%%%%%%%%%%%%%%%%%%%%%%%%%%%%%%%%%%%%%%%%%%%%%%%%%%%%%%%
%
%                                             T H E O R E M  
%                                       GENERALIZED ROUCHE - MATRICES
%
%%%%%%%%%%%%%%%%%%%%%%%%%%%%%%%%%%%%%%%%%%%%%%%%%%%%%%%%%%%%%%%%%%%%%%%%%%%%%%%%%%%%%%%%%%%%%%%%%%%%%%%%%%%%%%%%%%%%%%%%%%%%%%%%%%%%%%%%%%%%%%%%%%%%%%%%%%%%%%%%%%
%%%%%%%%%%%%%%%%%%%%%%%%%%%%%%%%%%%%%%%%%%%%%%%%%%%%%%%%%%%%%%%%%%%%%%%%%%%%%%%%%%%%%%%%%%%%%%%%%%%%%%%%%%%%%%%%%%%%%%%%%%%%%%%%%%%%%%%%%%%%%%%%%%%%%%%%%%%%%%%%%%
\begin{theorem}
\label{genRoucheMatrices}           
Let $A,B:\Omega \rightarrow \cmm$ be analytic matrix-valued functions, where $\Omega$ is an open
connected subset of $\complex$ 
and assume that $A(z)$ is nonsingular for all $z$ on the simple closed curve $\Gamma \subseteq \Omega$. 

If, for a subordinate matrix norm, $||A(z)^{-1}B(z)|| < 1$ for all $z \in \Gamma$,
then $\det(A+B)$ and $\det(A)$ have the same number of zeros inside $\Gamma$, counting multiplicities.
\end{theorem}
%%%%%%%%%%%%%%%%%%%%%%%%%%%%%%%%%%%%%%%%%%%%%%%%%%%%%%%%%%%%%%%%%%%%%%%%%%%%%%%%%%%%%%%%%%%%%%%%%%%%%%%%%%%%%%%%%%%%%%%%%%%%%%%%%%%%%%%%%%%%%%%%%%%%%%%%%%%%%%%%%%
%%%%%%%%%%%%%%%%%%%%%%%%%%%%%%%%%%%%%%%%%%%%%%%%%%%%%%%%%%%%%%%%%%%%%%%%%%%%%%%%%%%%%%%%%%%%%%%%%%%%%%%%%%%%%%%%%%%%%%%%%%%%%%%%%%%%%%%%%%%%%%%%%%%%%%%%%%%%%%%%%%
Theorem~\ref{genRoucheMatrices} is a convenient (although not the only) way to prove the following generalization to matrix polynomials of a 
result by Cauchy from~1829 (\cite{Cauchy}, \cite[Theorem (27,1), p.122]{Marden}). It can be found in~\cite{BiniNoferiniSharify}, \cite{HighamTisseur}, 
and \cite{Melman_MatPol}.
%%%%%%%%%%%%%%%%%%%%%%%%%%%%%%%%%%%%%%%%%%%%%%%%%%%%%%%%%%%%%%%%%%%%%%%%%%%%%%%%%%%%%%%%%%%%%%%%%%%%%%%%%%%%%%%%%%%%%%%%%%%%%%%%%%%%%%%%%%%%%%%%%%%%%%%%%%%%%%%%%%
%%%%%%%%%%%%%%%%%%%%%%%%%%%%%%%%%%%%%%%%%%%%%%%%%%%%%%%%%%%%%%%%%%%%%%%%%%%%%%%%%%%%%%%%%%%%%%%%%%%%%%%%%%%%%%%%%%%%%%%%%%%%%%%%%%%%%%%%%%%%%%%%%%%%%%%%%%%%%%%%%%
%
%                                             T H E O R E M  
%                                           GENERALIZED CAUCHY                 
%
%%%%%%%%%%%%%%%%%%%%%%%%%%%%%%%%%%%%%%%%%%%%%%%%%%%%%%%%%%%%%%%%%%%%%%%%%%%%%%%%%%%%%%%%%%%%%%%%%%%%%%%%%%%%%%%%%%%%%%%%%%%%%%%%%%%%%%%%%%%%%%%%%%%%%%%%%%%%%%%%%%
%%%%%%%%%%%%%%%%%%%%%%%%%%%%%%%%%%%%%%%%%%%%%%%%%%%%%%%%%%%%%%%%%%%%%%%%%%%%%%%%%%%%%%%%%%%%%%%%%%%%%%%%%%%%%%%%%%%%%%%%%%%%%%%%%%%%%%%%%%%%%%%%%%%%%%%%%%%%%%%%%%
\begin{theorem}
\label{genCauchy}             
The eigenvalues of the matrix polynomial $P(z) = A_{n}z^{n} + A_{n-1}z^{n-1} + \dots + A_{1}z + A_{0}$, with $A_{j} \in \cmm$ and $\An$ nonsingular,
are contained in the disk $|z| \leq \rho$, where $\rho$ is the unique positive root of
\bdis
||A^{-1}_{n}||^{-1}x^{n} - ||A_{n-1}||x^{n-1} - \dots - ||A_{1}||x - ||A_{0}|| = 0 \; ,
\edis
for any subordinate matrix norm.
\end{theorem}
%%%%%%%%%%%%%%%%%%%%%%%%%%%%%%%%%%%%%%%%%%%%%%%%%%%%%%%%%%%%%%%%%%%%%%%%%%%%%%%%%%%%%%%%%%%%%%%%%%%%%%%%%%%%%%%%%%%%%%%%%%%%%%%%%%%%%%%%%%%%%%%%%%%%%%%%%%%%%%%%%%
%%%%%%%%%%%%%%%%%%%%%%%%%%%%%%%%%%%%%%%%%%%%%%%%%%%%%%%%%%%%%%%%%%%%%%%%%%%%%%%%%%%%%%%%%%%%%%%%%%%%%%%%%%%%%%%%%%%%%%%%%%%%%%%%%%%%%%%%%%%%%%%%%%%%%%%%%%%%%%%%%%
This theorem leads to the following definition.
%%%%%%%%%%%%%%%%%%%%%%%%%%%%%%%%%%%%%%%%%%%%%%%%%%%%%%%%%%%%%%%%%%%%%%%%%%%%%%%%%%%%%%%%%%%%%%%%%%%%%%%%%%%%%%%%%%%%%%%%%%%%%%%%%%%%%%%%%%%%%%%%%%%%%%%%%%%%%%%%%%
%%%%%%%%%%%%%%%%%%%%%%%%%%%%%%%%%%%%%%%%%%%%%%%%%%%%%%%%%%%%%%%%%%%%%%%%%%%%%%%%%%%%%%%%%%%%%%%%%%%%%%%%%%%%%%%%%%%%%%%%%%%%%%%%%%%%%%%%%%%%%%%%%%%%%%%%%%%%%%%%%%
%
%                                             D E F I N I T I O N
%                                                CAUCHY RADIUS                
%
%%%%%%%%%%%%%%%%%%%%%%%%%%%%%%%%%%%%%%%%%%%%%%%%%%%%%%%%%%%%%%%%%%%%%%%%%%%%%%%%%%%%%%%%%%%%%%%%%%%%%%%%%%%%%%%%%%%%%%%%%%%%%%%%%%%%%%%%%%%%%%%%%%%%%%%%%%%%%%%%%%
%%%%%%%%%%%%%%%%%%%%%%%%%%%%%%%%%%%%%%%%%%%%%%%%%%%%%%%%%%%%%%%%%%%%%%%%%%%%%%%%%%%%%%%%%%%%%%%%%%%%%%%%%%%%%%%%%%%%%%%%%%%%%%%%%%%%%%%%%%%%%%%%%%%%%%%%%%%%%%%%%%
\begin{definition}[Cauchy radius]
The quantity $\rho$ in the previous theorem is called the Cauchy radius of the matrix polynomial $P$. It depends on the matrix norm used
in the theorem.
\end{definition}
%%%%%%%%%%%%%%%%%%%%%%%%%%%%%%%%%%%%%%%%%%%%%%%%%%%%%%%%%%%%%%%%%%%%%%%%%%%%%%%%%%%%%%%%%%%%%%%%%%%%%%%%%%%%%%%%%%%%%%%%%%%%%%%%%%%%%%%%%%%%%%%%%%%%%%%%%%%%%%%%%%
%%%%%%%%%%%%%%%%%%%%%%%%%%%%%%%%%%%%%%%%%%%%%%%%%%%%%%%%%%%%%%%%%%%%%%%%%%%%%%%%%%%%%%%%%%%%%%%%%%%%%%%%%%%%%%%%%%%%%%%%%%%%%%%%%%%%%%%%%%%%%%%%%%%%%%%%%%%%%%%%%%
Theorem~\ref{genCauchy} will be the reference inclusion region to which we will compare our results, since it generally appears to be among the best bounds 
attainable for matrix polynomials expressed in the standard power basis, judging from the extensive results in~\cite{HighamTisseur}, where a large number 
of such eigenvalue bounds were compared.

Throughout the paper, we will use $I$ for the identity matrix without specifying its size, which is usually clear from the context. On those occasions
where it is not, a $k \times k$ identity matrix will be denoted by $I_{k}$. 

The paper is organized as follows. In Section~\ref{mainresult} we derive an inclusion region for the eigenvalues of a matrix polynomial expressed in a 
generalized basis, which is then applied to several engineering problems in Section~\ref{examples}.

%%%%%%%%%%%%%%%%%%%%%%%%%%%%%%%%%%%%%%%%%%%%%%%%%%%%%%%%%%%%%%%%%%%%%%%%%%%%%%%%%%%%%%%%%%%%%%%%%%%%%%%%%%%%%%%%%%%%%%%%%%%%%%%%%%%%%%%%%%%%%%%%%%%%%%%%%%%%%%%%%%
%%%%%%%%%%%%%%%%%%%%%%%%%%%%%%%%%%%%%%%%%%%%%%%%%%%%%%%%%%%%%%%%%%%%%%%%%%%%%%%%%%%%%%%%%%%%%%%%%%%%%%%%%%%%%%%%%%%%%%%%%%%%%%%%%%%%%%%%%%%%%%%%%%%%%%%%%%%%%%%%%%
%%%%%%%%%%%%%%%%%%%%%%%%%%%%%%%%%%%%%%%%%%%%%%%%%%%%%%%%%%%%%%%%%%%%%%%%%%%%%%%%%%%%%%%%%%%%%%%%%%%%%%%%%%%%%%%%%%%%%%%%%%%%%%%%%%%%%%%%%%%%%%%%%%%%%%%%%%%%%%%%%%
%
%
%              SECTION 2  -  MAIN RESULT                                               
%
%
%%%%%%%%%%%%%%%%%%%%%%%%%%%%%%%%%%%%%%%%%%%%%%%%%%%%%%%%%%%%%%%%%%%%%%%%%%%%%%%%%%%%%%%%%%%%%%%%%%%%%%%%%%%%%%%%%%%%%%%%%%%%%%%%%%%%%%%%%%%%%%%%%%%%%%%%%%%%%%%%%%
%%%%%%%%%%%%%%%%%%%%%%%%%%%%%%%%%%%%%%%%%%%%%%%%%%%%%%%%%%%%%%%%%%%%%%%%%%%%%%%%%%%%%%%%%%%%%%%%%%%%%%%%%%%%%%%%%%%%%%%%%%%%%%%%%%%%%%%%%%%%%%%%%%%%%%%%%%%%%%%%%%
%%%%%%%%%%%%%%%%%%%%%%%%%%%%%%%%%%%%%%%%%%%%%%%%%%%%%%%%%%%%%%%%%%%%%%%%%%%%%%%%%%%%%%%%%%%%%%%%%%%%%%%%%%%%%%%%%%%%%%%%%%%%%%%%%%%%%%%%%%%%%%%%%%%%%%%%%%%%%%%%%%

\section{Main result}   
\label{mainresult}

%%%%%%%%%%%%%%%%%%%%%%%%%%%%%%%%%%%%%%%%%%%%%%%%%%%%%%%%%%%%%%%%%%%%%%%%%%%%%%%%%%%%%%%%%%%%%%%%%%%%%%%%%%%%%%%%%%%%%%%%%%%%%%%%%%%%%%%%%%%%%%%%%%%%%%%%%%%%%%%%%%
%%%%%%%%%%%%%%%%%%%%%%%%%%%%%%%%%%%%%%%%%%%%%%%%%%%%%%%%%%%%%%%%%%%%%%%%%%%%%%%%%%%%%%%%%%%%%%%%%%%%%%%%%%%%%%%%%%%%%%%%%%%%%%%%%%%%%%%%%%%%%%%%%%%%%%%%%%%%%%%%%%
%
%
%                                                          THEOREM - MAIN                                                  
%
%
%%%%%%%%%%%%%%%%%%%%%%%%%%%%%%%%%%%%%%%%%%%%%%%%%%%%%%%%%%%%%%%%%%%%%%%%%%%%%%%%%%%%%%%%%%%%%%%%%%%%%%%%%%%%%%%%%%%%%%%%%%%%%%%%%%%%%%%%%%%%%%%%%%%%%%%%%%%%%%%%%%
%%%%%%%%%%%%%%%%%%%%%%%%%%%%%%%%%%%%%%%%%%%%%%%%%%%%%%%%%%%%%%%%%%%%%%%%%%%%%%%%%%%%%%%%%%%%%%%%%%%%%%%%%%%%%%%%%%%%%%%%%%%%%%%%%%%%%%%%%%%%%%%%%%%%%%%%%%%%%%%%%%
\begin{theorem}
\label{maintheorem}
Let $\{q_{j}\}_{j=0}^{n}$ be a scalar polynomial basis, where $q_{j}$ is a polynomial of degree $j$ for $j=1,\dots,n$, and
denote by $r_{ij}$ the $i$th zero of $q_{j}$.
If for every $j \geq 0$ there exist
nonnegative numbers $\alpha_{1}^{(j)},\dots,\alpha_{j}^{(j)}$ so that $\sum_{i=1}^{j} \alpha_{i}^{(j)} \leq \gamma$, with $\gamma > 0$, 
and, for $z \neq r_{ij}$,
\beq
\label{qineq0}
\left | \dfrac{q_{j-1}(z)}{q_{j}(z)} \right | 
\leq \sum_{i=1}^{j} \dfrac{\alpha_{i}^{(j)}}{|z-r_{ij}|} \; , 
\eeq
then the eigenvalues of the matrix polynomial 
\bdis
P(z) = \An  \qn(z) + \Anmo \qnmo(z) + \dots \Aone \qone(z) + \Azero \qzero(z) \; ,
\edis
with $\Aj \in \complex^{m\times m}$ and $\An$ nonsingular, are contained in the union of the at most $n(n+1)/2$ disks
\bdis
\mathcal{R} = \stackrel{n}{\underset{i \leq j}{\underset{i,j = 1}{\bigcup}}}  \lcb z \in \complex : |z-r_{ij}| \leq \gamma \rho \rcb \; ,
\edis
where $\rho$ is the \emph{Cauchy radius} of $\, \sum_{j=0}^{n} \Aj z^{j}$.
Moreover, if the region $\mathcal{R}$ is composed of disjoint components, then each component contains $m$ times as many eigenvalues of $P$ as 
it contains zeros of $\qn$.
\end{theorem}
%%%%%%%%%%%%%%%%%%%%%%%%%%%%%%%%%%%%%%%%%%%%%%%%%%%%%%%%%%%%%%%%%%%%%%%%%%%%%%%%%%%%%%%%%%%%%%%%%%%%%%%%%%%%%%%%%%%%%%%%%%%%%%%%%%%%%%%%%%%%%%%%%%%%%%%%%%%%%%%%%%
\prf  All matrix norms are considered to be subordinate. 
If $z$ is an eigenvalue of $P$ such that $z \neq r_{ij}$, then 
\bdis
\text{det} \bigl (  \An \qn(z) + \dots + \Aone \qone(z) + \Azero \qzero(z) \bigr )  = 0 
\edis
implies that 
\bdis
\text{det} \biggl ( I + \bigl ( \An \qn(z) \bigr )^{-1} \bigl (  \Anmo \qnmo(z) + \dots + \Aone \qone(z) + \Azero \qzero(z) \bigr ) \biggr ) = 0 \; ,
\edis     
which is only possible (\cite[p.351]{HJ}) if 
\beq
\label{detineq1}
\lnorm \bigl ( \An \qn(z) \bigr )^{-1} \bigl ( \Anmo \qnmo(z) + \dots + \Aone \qone(z) + \Azero \qzero(z) \bigr ) \rnorm \geq 1 \; .
\eeq
Since $\|A^{-1}\|\|B\| \geq \|A^{-1}B\|$, inequality~(\ref{detineq1}) implies that 
\beq
\label{detineq2}
\lnorm \Anmo \qnmo(z) + \dots + \Aone \qone(z) + \Azero \qzero(z) \rnorm \geq \lnorm \An^{-1} \rnorm^{-1} |\qn(z)| \; ,
\eeq
so that, with $\|A\| + \|B\| \geq \|A + B \|$, inequality~(\ref{detineq2}) yields
\beq
\label{detineq3}                    
\lnorm \Anmo \rnorm |\qnmo(z)| + \dots + \lnorm \Aone \rnorm |\qone(z)| + \lnorm \Azero \rnorm |\qzero(z)| \geq \lnorm \An^{-1} \rnorm ^{-1} |q_{n}(z)| \; .
\eeq

To express the left-hand side of~(\ref{detineq3}) in a more useful way, we define for each $j$ ($j=1,2,\dots,n$):
\bdis
d_{j}(z) = \underset{i \leq k}{\underset{1 \leq i,k \leq j}{\text{min}}} |z-r_{ik}|  \; ,
\edis
namely, the distance of a point $z$ to set of all the zeros of $\,q_{1},\dots,q_{j}$. Clearly, 
$d_{1}(z) \geq d_{2}(z) \geq \dots \geq d_{n}(z)$. 
Inequality~(\ref{qineq0}) then implies that
\beq
\label{qineq1}
\left | \dfrac{q_{j-1}(z)}{q_{j}(z)} \right | 
\leq \sum_{i=1}^{j} \dfrac{\alpha_{i}^{(j)}}{|z-r_{ij}|} 
\leq \dfrac{\sum_{i=1}^{j} \alpha_{i}^{(j)}}{\underset{1 \leq i \leq j}{\text{min}} |z-r_{ij}|} 
\leq \dfrac{\gamma}{d_{j}(z)} 
\leq \dfrac{\gamma}{d_{n}(z)}  
\qquad \text{($j=1,\dots,n$).}
\eeq
Repeated application of~(\ref{qineq1}) yields
\beq
\label{qineq2}
\left | \dfrac{q_{j}(z)}{q_{n}(z)} \right | 
= \left | \dfrac{q_{j}(z)}{q_{j+1}(z)} \right | 
\dots
\left | \dfrac{q_{n-1}(z)}{q_{n}(z)} \right | 
\leq \lb \dfrac{\gamma}{d_{n}(z)} \rb^{n-j} \; \cdot
\eeq

Dividing~(\ref{detineq3}) by $|\qn(z)|$ and majorizing its left-hand side in terms of $d_{n}(z)$ using~(\ref{qineq2}) yields
\begin{eqnarray}
& & \hskip -1.50cm \lnorm \Anmo \rnorm \left | \dfrac{\qnmo(z)}{\qn(z)} \right | + \dots + \lnorm \Aone \rnorm \left | \dfrac{\qone(z)}{\qn(z)} \right | 
                   + \lnorm \Azero \rnorm \left | \dfrac{\qzero(z)}{\qn(z)} \right | \nonumber \\
& & \hskip 2cm \leq \lnorm \Anmo \rnorm \lb \dfrac{\gamma}{d_{n}(z)} \rb + \dots + \lnorm \Aone \rnorm \lb \dfrac{\gamma}{d_{n}(z)} \rb^{n-1} 
                    + \lnorm \Azero \rnorm \lb \dfrac{\gamma}{d_{n}(z)} \rb ^{n}  \nonumber \\
& & \hskip 2cm \leq \lb \, \lnorm \Anmo \rnorm \lb \dfrac{d_{n}(z)}{\gamma} \rb ^{n-1} + \dots + \lnorm \Aone \rnorm \lb \dfrac{d_{n}(z)}{\gamma} \rb 
                + \lnorm \Azero \rnorm \, \rb \lb \dfrac{\gamma}{d_{n}(z)} \rb ^{n} \; \cdot \label{polynomineq1}
\end{eqnarray} 
Combining~(\ref{polynomineq1}) with~(\ref{detineq3}), we obtain that if $z$ is an eigenvalue of $P$, then 
\bdis
\lnorm \Anmo \rnorm \lb \dfrac{d_{n}(z)}{\gamma} \rb ^{n-1} + \dots + \lnorm \Aone \rnorm \lb \dfrac{d_{n}(z)}{\gamma} \rb + \lnorm \Azero \rnorm 
\geq
\lnorm \An^{-1} \rnorm ^{-1} \lb \dfrac{d_{n}(z)}{\gamma} \rb ^{n} \; \cdot 
\edis
As can be seen from Theorem~\ref{genCauchy}, this means that $d_{n}(z)/\gamma \leq \rho$, where $\rho$ is the Cauchy radius of $\sum_{j=0}^{n}\Aj z^{j}$, or 
\bdis
\underset{i \leq k}{\underset{1 \leq i,k \leq n}{\text{min}}} |z-r_{ik}|  \leq \gamma \rho \; ,
\edis
i.e., $z$ must lie in the union $\mathcal{R}$ of disks centered at the zeros of $\qone, \dots, \qn$ with radius $\gamma \rho$.
The number of those disks is at most $\sum_{j=1}^{n} j = n(n+1)/2$,
since basis polynomials can have common zeros.

If $\mathcal{R}$ is composed of disjoint subregions (each subregion necessarily a union of disks), then the boundary $\Gamma$ of such a subregion
is a simple closed curve on which $\An \qn(z)$ is nonsingular and $d_{n}(z) = \gamma \rho$. 
Now consider the collection of points $z$ for which $d_{n}(z) = \gamma \rho + \varepsilon$ for $\varepsilon > 0$. It is the boundary of 
the union $\mathcal{R}_{1}$ of disks with the same centers as those that determine $\mathcal{R}$, but with a larger radius. Clearly, this boundary
does not contain any of the centers $r_{ij}$. If $\mathcal{R}$ consists of disjoint subregions, then we can choose $\varepsilon$ small enough so that 
$\mathcal{R}_{1}$ does as well. One of those will necessarily enclose $\Gamma$, and we define $\Gamma_{1}$ as its boundary.
It is a simple closed curve on which $\An \qn(z)$ is nonsingular and $d_{n}(z) = \gamma \rho + \varepsilon$.  

Using the same arguments as in~(\ref{detineq2}),~(\ref{detineq3}), and~(\ref{polynomineq1}), one sees that, for $z \neq r_{ij}$, the inequality 
\beq        
\label{Roucheineq1}
\lnorm \bigl ( \An q_{n}(z) \bigr ) ^{-1} \bigl ( \Anmo \qnmo(z) + \dots + \Aone \qone(z) + \Azero \qzero(z) \bigr ) \rnorm < 1 
\eeq
will certainly 
be satisfied when 
\beq
\label{Roucheineq2}
\lnorm \Anmo \rnorm \lb \dfrac{d_{n}(z)}{\gamma} \rb ^{n-1} + \dots + \lnorm \Aone \rnorm \lb \dfrac{d_{n}(z)}{\gamma} \rb + \lnorm \Azero \rnorm 
<
\lnorm \An^{-1} \rnorm ^{-1} \lb \dfrac{d_{n}(z)}{\gamma} \rb ^{n}  \; .         
\eeq
Since for any $z \in \Gamma_{1}$ we have $d_{n}(z)/\gamma = \rho + \varepsilon/\gamma > \rho$, Theorem~\ref{genCauchy} implies that~(\ref{Roucheineq2})
is satisfied on~$\Gamma_{1}$. This, in turn, implies that~(\ref{Roucheineq1}) is satisfied, from which we obtain with Theorem~\ref{genRoucheMatrices} 
that $P$ and $\An \qn$ have 
the same number of eigenvalues in the open region enclosed by $\Gamma_{1}$. Since $\Gamma_{1}$ encloses $\Gamma$, we conclude, by letting 
$\varepsilon \rightarrow 0^{+}$, that the closed subregion of $\mathcal{R}$ bounded by $\Gamma$ contains a number of eigenvalues of $P$
equal to the number of eigenvalues of $\An \qn$ that it contains.
Because $\text{det} (\An  \qn(w) )~=~0 \Longleftrightarrow \text{det}(\An) \, \qn^{m}(w)~=~0$, and $\text{det}(\An)~\neq~0$, 
this number is $m$ times the number of zeros of $\qn$ in the closed region. \qed
%%%%%%%%%%%%%%%%%%%%%%%%%%%%%%%%%%%%%%%%%%%%%%%%%%%%%%%%%%%%%%%%%%%%%%%%%%%%%%%%%%%%%%%%%%%%%%%%%%%%%%%%%%%%%%%%%%%%%%%%%%%%%%%%%%%%%%%%%%%%%%%%%%%%%%%%%%%%%%%%%%
%%%%%%%%%%%%%%%%%%%%%%%%%%%%%%%%%%%%%%%%%%%%%%%%%%%%%%%%%%%%%%%%%%%%%%%%%%%%%%%%%%%%%%%%%%%%%%%%%%%%%%%%%%%%%%%%%%%%%%%%%%%%%%%%%%%%%%%%%%%%%%%%%%%%%%%%%%%%%%%%%%
\vskip 0.5cm

In the special case where $P$ is a scalar polynomial and the polynomials $\qj$ form a weakly interlacing system, i.e., polynomials with weakly 
interlacing real zeros, Theorem~\ref{maintheorem} essentially reduces to Theorem~8.4.6 in~\cite{RS}. In this case, Lemma~8.4.5 in~\cite{RS} shows that
condition~(\ref{qineq0}) is statisfied with $\gamma=1$, the zeros $r_{ij}$ all lie on the real axis, and the region derived in Theorem~8.4.6 in~\cite{RS} 
is the convex hull of the one in Theorem~\ref{maintheorem}. Weakly interlacing bases include all classical orthogonal bases: Hermite, Legendre, Chebyshev, etc.

Another special case is obtained by choosing the Newton basis with complex nodes $\{a_{j}\}$, $j=1,2,\dots$, for which the basis polynomials are defined by
$q_{0}(z)=1$ and $q_{j}(z)=(z-a_{j})q_{j-1}$. This means that the zeros of $q_{j}$ are $a_{1},a_{2},\dots,a_{j}$.
Since
\bdis
\dfrac{q_{j-1}(z)}{q_{j}(z)} = \dfrac{1}{z-a_{j}} \; ,
\edis
condition~(\ref{qineq0}) is satisfied with $\alpha_{1}^{(j)}=\alpha_{2}^{(j}=\dots=\alpha_{j-1}^{(j)}=0$, $\alpha_{j}^{(j)}=1$, and $\gamma=1$.
For scalar polynomials, this is Theorem~8.6.3 in~\cite{RS}.
However, Theorem~\ref{maintheorem} allows for more general bases where the zeros of different basis polynomials do not need to be related.

Obviously, changing the basis does not universally improve results for all problems. As with other localization results, some problems lend themselves
better to certain bounds than others.
%but problems for which it can be useful are relatively easy to characterize. 
However, the degree of matrix polynomials appearing in engineering applications tends to be low; many of them are quadratic. They are easily expressed 
in a different basis, while the computation of bounds 
is several orders of magnitude less onerous than the computation of the actual eigenvalues, so that not much is lost by trying a different basis.
To illustrate how Theorem~\ref{maintheorem} can improve classical bounds, we turn to the literature on quadratic eigenvalue problems with their many
applications in engineering.

%%%%%%%%%%%%%%%%%%%%%%%%%%%%%%%%%%%%%%%%%%%%%%%%%%%%%%%%%%%%%%%%%%%%%%%%%%%%%%%%%%%%%%%%%%%%%%%%%%%%%%%%%%%%%%%%%%%%%%%%%%%%%%%%%%%%%%%%%%%%%%%%%%%%%%%%%%%%%%%%%%
%%%%%%%%%%%%%%%%%%%%%%%%%%%%%%%%%%%%%%%%%%%%%%%%%%%%%%%%%%%%%%%%%%%%%%%%%%%%%%%%%%%%%%%%%%%%%%%%%%%%%%%%%%%%%%%%%%%%%%%%%%%%%%%%%%%%%%%%%%%%%%%%%%%%%%%%%%%%%%%%%%
%%%%%%%%%%%%%%%%%%%%%%%%%%%%%%%%%%%%%%%%%%%%%%%%%%%%%%%%%%%%%%%%%%%%%%%%%%%%%%%%%%%%%%%%%%%%%%%%%%%%%%%%%%%%%%%%%%%%%%%%%%%%%%%%%%%%%%%%%%%%%%%%%%%%%%%%%%%%%%%%%%
%
%
%              SECTION 3  -  EXAMPLES                                                  
%
%
%%%%%%%%%%%%%%%%%%%%%%%%%%%%%%%%%%%%%%%%%%%%%%%%%%%%%%%%%%%%%%%%%%%%%%%%%%%%%%%%%%%%%%%%%%%%%%%%%%%%%%%%%%%%%%%%%%%%%%%%%%%%%%%%%%%%%%%%%%%%%%%%%%%%%%%%%%%%%%%%%%
%%%%%%%%%%%%%%%%%%%%%%%%%%%%%%%%%%%%%%%%%%%%%%%%%%%%%%%%%%%%%%%%%%%%%%%%%%%%%%%%%%%%%%%%%%%%%%%%%%%%%%%%%%%%%%%%%%%%%%%%%%%%%%%%%%%%%%%%%%%%%%%%%%%%%%%%%%%%%%%%%%
%%%%%%%%%%%%%%%%%%%%%%%%%%%%%%%%%%%%%%%%%%%%%%%%%%%%%%%%%%%%%%%%%%%%%%%%%%%%%%%%%%%%%%%%%%%%%%%%%%%%%%%%%%%%%%%%%%%%%%%%%%%%%%%%%%%%%%%%%%%%%%%%%%%%%%%%%%%%%%%%%%

\section{Examples}
\label{examples}

We establish a few preliminary results concerning monic quadratic matrix polynomials before applying them to numerical examples.
Since it is easy to compute, we choose the $1$-norm throughout this section.

%%%%%%%%%%%%%%%%%%%%%%%%%%%%%%%%%%%%%%%%%%%%%%%%%%%%%%%%%%%%%%%%%%%%%%%%%%%%%%%%%%%%%%%%%%%%%%%%%%%%%%%%%%%%%%%%%%%%%%%%%%%%%%%%%%%%%%%%%%%%%%%%%%%%%%%%%%%%%%%%%%
%
%
%              SUBSECTION 3.1  -  QUADRATIC MATRIX POLYNOMIALS                                                  
%
%
%%%%%%%%%%%%%%%%%%%%%%%%%%%%%%%%%%%%%%%%%%%%%%%%%%%%%%%%%%%%%%%%%%%%%%%%%%%%%%%%%%%%%%%%%%%%%%%%%%%%%%%%%%%%%%%%%%%%%%%%%%%%%%%%%%%%%%%%%%%%%%%%%%%%%%%%%%%%%%%%%%

\subsection {Quadratic matrix polynomials}   
\label{quadraticproblems}

We consider the monic quadratic matrix polynomial $P(z) = Iz^{2} + \Aone z + \Azero$, expressed in the standard power basis, and we define
the Newton basis $\mathcal{N} = \{ \fzero, \fone, \ftwo \}$ and the more general basis $\mathcal{B} = \{ \qzero, \qone, \qtwo \}$, respectively, by
\bdis
\begin{array}{ccc}
         \left \{ \begin{array}{l}
                  \fzeroz = 1               \\
                  \fonez = z-a              \\
                  \ftwoz = (z-a)(z-b)       \\
                  \end{array}
         \right .
                                                          & \text{and} &  
                                                          \hskip 0.5cm     \left \{ \begin{array}{l}
                                                                                    \qzeroz = 1               \\
                                                                                    \qonez = z-a              \\
                                                                                    \qtwoz = (z-b)(z-c) \; ,  \\ 
                                                                                    \end{array}
                                                                                    \right .
\end{array}
\edis
where for the basis $\mathcal{B}$, either $a=b=c$ in which case it becomes a Newton basis, or $b \neq c$.
The choice of the nodes $a$, $b$, and $c$, which are generally different for different bases, is taylored to the particular matrix polynomial.
The power basis is easily expressed in terms of $\mathcal{N}$ and $\mathcal{B}$: 
\bdis
\hskip -0.5cm \begin{array}{ccc}
         \left \{ \begin{array}{l}
                  1 = \fzeroz                \\
                  z = \fonez + a \fzeroz     \\
                  z^{2} = \ftwoz + (a+b) \fonez + a^{2} \fzeroz      \\
                  \end{array}
         \right .
                                                          & \hskip -0.25cm \text{and} &  
                                                          \hskip 0.15cm     \left \{ \begin{array}{l}
                                                                                    1 = \qzeroz               \\
                                                                                    z = \qonez + a \qzeroz    \\
                                                                                    z^{2} = \qtwoz + (b+c)\qonez + (a(b+c) -bc)\qzeroz \; .     \\
                                                                                    \end{array}
                                                                                    \right .
\end{array}
\edis
The quadratic $P$ in the bases $\mathcal{N}$ and $\mathcal{B}$ then becomes, respectively,
\begin{eqnarray}
P(z) & = & I\ftwoz + \bigl ( \Aone + (a+b)I \bigr ) \fonez + \bigl ( a\Aone + \Azero + a^{2} I \bigr ) \fzeroz \; , \nonumber  \\
P(z) & = & I\qtwoz + \bigl ( \Aone + (b+c)I \bigr ) \qonez + \bigl ( a\Aone + \Azero + (a(b+c) -bc)I \bigr ) \qzeroz \; .  \nonumber \\
\label{pinbases} 
\end{eqnarray}
Let us verify condition~(\ref{qineq0}) in Theorem~\ref{maintheorem} for these bases. For $\mathcal{N}$, this was already done in the remarks following 
that theorem. For the basis $\mathcal{B}$ with $b \neq c$, we obtain
\bdis
\dfrac{\qzeroz}{\qonez} = \dfrac{1}{z-a}
\;\; \text{and} \;\;
\dfrac{\qonez}{\qtwoz} = \dfrac{(a-b)/(c-b)}{z-b} + \dfrac{(c-a)/(c-b)}{z-c} \; , 
\edis
so that~(\ref{qineq0}) is satisfied with 
\bdis
\alpha_{1}^{(1)}=1 \; , \; \alpha_{1}^{(2)} = \left | \dfrac{a-b}{c-b} \right | \; , \; \alpha_{2}^{(2)} = \left | \dfrac{a-c}{c-b} \right | \; , \; \text{and} \; 
\gamma = \left | \dfrac{a-b}{c-b} \right | + \left | \dfrac{a-c}{c-b} \right | \; \cdot
\edis
We note that $\gamma \geq 1$.

The nodes determining the bases $\mathcal{N}$ and $\mathcal{B}$ should be chosen so as to make the norms of the coefficient matrices as small as possible,
since this will make the radii of the disks in the inclusion region smaller. To do this for the numerical examples below, we will use the observation that 
for real numbers $\{\beta_{j}\}_{j=1}^{n}$, ordered in increasing order, the solution of the minimization problem 
\beq
\label{minproblem}
\min_{x} \max_{1 \leq j \leq n} |\beta_{j} -x | 
\eeq
is obtained for $x^{*}=(\beta_{1}+\beta_{n})/2$. This implies that if the numbers $\beta_{j}$ are the diagonal of a diagonal matrix $M$, then 
$\|M-x^{*}I\|_{1} \leq \|M\|_{1}$. When the matrix $M$ is not diagonal, but strongly diagonally dominant, then we expect this inequality to 
still be true in most cases. When the numbers $\beta_{j}$ are complex and the minimization in $x$ is to be carried out over the complex plane,
then, to keep matters simple, we will carry out the minimization separately for the real and complex parts.  

%%%%%%%%%%%%%%%%%%%%%%%%%%%%%%%%%%%%%%%%%%%%%%%%%%%%%%%%%%%%%%%%%%%%%%%%%%%%%%%%%%%%%%%%%%%%%%%%%%%%%%%%%%%%%%%%%%%%%%%%%%%%%%%%%%%%%%%%%%%%%%%%%%%%%%%%%%%%%%%%%%

%%%%%%%%%%%%%%%%%%%%%%%%%%%%%%%%%%%%%%%%%%%%%%%%%%%%%%%%%%%%%%%%%%%%%%%%%%%%%%%%%%%%%%%%%%%%%%%%%%%%%%%%%%%%%%%%%%%%%%%%%%%%%%%%%%%%%%%%%%%%%%%%%%%%%%%%%%%%%%%%%%
%
%
%              SUBSECTION 3.2  -  NUMERICAL EXAMPLES                                                            
%
%
%%%%%%%%%%%%%%%%%%%%%%%%%%%%%%%%%%%%%%%%%%%%%%%%%%%%%%%%%%%%%%%%%%%%%%%%%%%%%%%%%%%%%%%%%%%%%%%%%%%%%%%%%%%%%%%%%%%%%%%%%%%%%%%%%%%%%%%%%%%%%%%%%%%%%%%%%%%%%%%%%%

\subsection {Numerical examples}                    
\label{numericalexamples} 
%%%%%%%%%%%%%%%%%%%%%%%%%%%%%%%%%%%%%%%%%%%%%%%%%%%%%%%%%%%%%%%%%%%%%%%%%%%%%%%%%%%%%%%%%%%%%%%%%%%%%%%%%%%%%%%%%%%%%%%%%%%%%%%%%%%%%%%%%%%%%%%%%%%%%%%%%%%%%%%%%%
%
%
%                                                  EXAMPLE 1                                                               
%
%
%%%%%%%%%%%%%%%%%%%%%%%%%%%%%%%%%%%%%%%%%%%%%%%%%%%%%%%%%%%%%%%%%%%%%%%%%%%%%%%%%%%%%%%%%%%%%%%%%%%%%%%%%%%%%%%%%%%%%%%%%%%%%%%%%%%%%%%%%%%%%%%%%%%%%%%%%%%%%%%%%%
{\bf Example 1.} We consider the connected damped mass-spring system in~\cite[p.259]{TisseurMeerbergen}. 
Its vibration is described by a second-order differential
equation of the form $A_{2}y''(t)+\Aone y'(t)+\Azero y(t) = f(t)$, where $A_{2}$, $\Aone$, and $\Azero$ are $m \times m$ matrices and $y(t)$ is an $m$-vector. 
The solution of the differential equation can be expressed in terms of the eigenvalues and eigenvectors of the quadratic eigenvalue 
problem $(A_{2}z^{2}+\Aone z+\Azero)v=0$.
Here, the mass matrix $A_{2}$ is diagonal, and the damping and stiffness matrices $\Aone$ and $\Azero$, respectively, are symmetric tridiagonal. 
In~\cite{TisseurMeerbergen}, $A_{2}=I$, $\Aone = \tau \, \text{tridiag}(-1,3,-1)$, $\Azero~=~\kappa \, \text{tridiag}(-1,3,-1)$, $\tau, \kappa \in \real$, 
and $m=50$. 
We will compare the standard power basis $\{1,z,z^{2}\}$ with the bases $\mathcal{N}$ and $\mathcal{B}$ from Subsection~\ref{quadraticproblems}.

We now determine the nodes $a$, $b$, and $c$ defining those bases, and start 
with the Newton basis $\mathcal{N}$, where from~(\ref{pinbases}) the coefficients of $\fone$ and $\fzero$ are given, respectively, by
\begin{eqnarray*}
\Aone + (a+b)I & = & \text{tridiag}(-\tau,3\tau+a+b,-\tau) \; ,   \\
a\Aone + \Azero + a^{2} I & = & \text{tridiag}(-\tau a -\kappa,a^{2}+3\tau a +3\kappa,-\tau a-\kappa) \; .
\end{eqnarray*}
In choosing $a$ and $b$, we aim to make the $1$-norm of the coefficients as small as possible. Without attempting an elaborate optimization, we will choose
$a$ and $b$ such as to make the diagonals of the matrix coefficients zero, i.e., $a^{2}+3\tau a +3\kappa=0$ and $a+b=-3\tau$. This means that $a$ and $b$
are the two zeros of $a^{2}+3\tau a +3\kappa$, and we choose $a$ as the zero for which $|a\tau + \kappa|$ is smaller.

For the basis $\mathcal{B}$, we have more flexibility since we now have three nodes $a$, $b$, and $c$. Here, from~(\ref{pinbases}), 
the coefficients of $\qone$ and $\qzero$
are given, respectively, by
\begin{eqnarray*}
\Aone + (b+c)I & = & \text{tridiag}(-\tau,3\tau+b+c,-\tau) \; ,   \\
a\Aone + \Azero + (a(b+c) -bc)I & = & \text{tridiag}(-\tau a-\kappa,a(b+c)-bc+3\tau a  +3\kappa,-\tau a-\kappa) \; .
\end{eqnarray*}
Arguing similarly as before, we choose the nodes such that $b+c=-3\tau$, $bc=3\kappa$ and $a=-\kappa/\tau$. This makes the diagonal of the coefficient
matrix of $\qone$ zero, while making the coefficient matrix of $\qzero$ vanish. The nodes $b$ and $c$ are the same as the nodes $a$ and $b$ we found for
the Newton basis since they are the zeros of the same quadratic. If $b=c$, then we set $a=b=c$, reverting to a Newton basis. Potentially better results could
obtained than for the Newton basis when $b \neq c$ although there is a price to pay in the form of a larger value for $\gamma$. 
It is therefore not a priori clear which basis is preferable. Fortunately, it is a simple matter to compute the $1$-norm, so that both bases can easily be 
compared. 

\begin{figure}[H]
\begin{center}
\raisebox{0ex}{\includegraphics[width=0.35\linewidth]{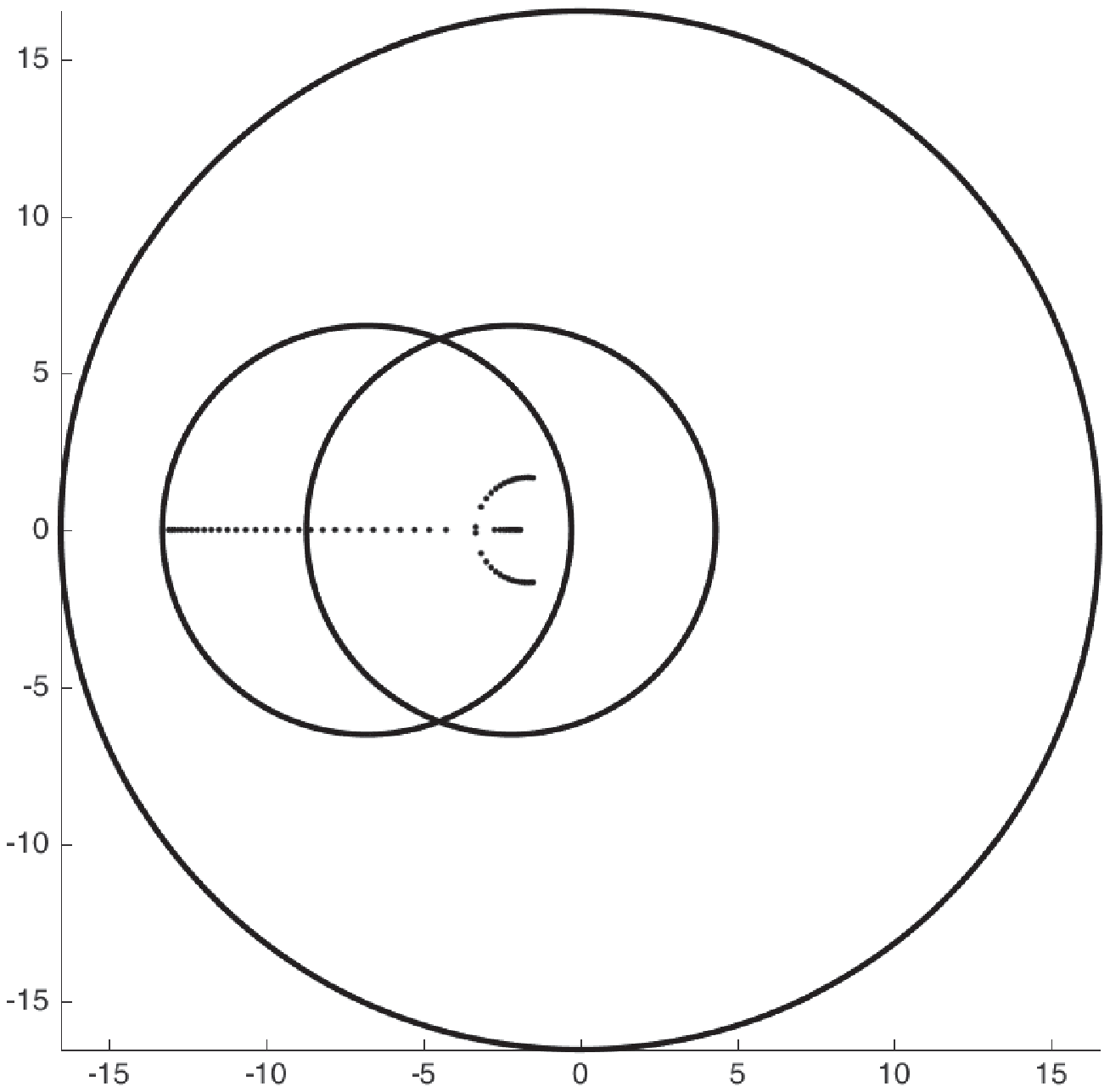}}
\hskip 2cm
\raisebox{0ex}{\includegraphics[width=0.35\linewidth]{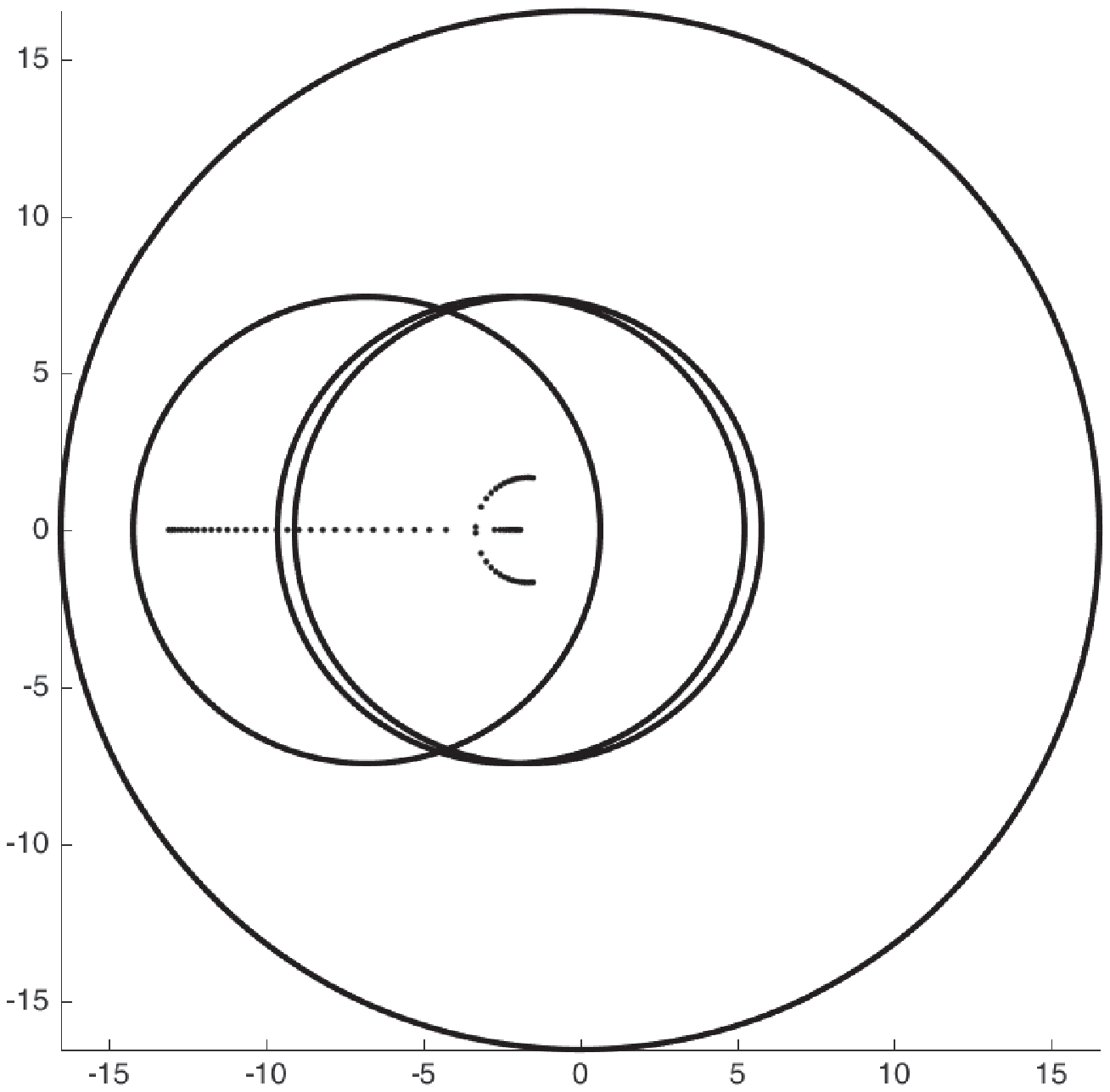}}
\caption{Inclusion regions for Example~1 with $\tau=3$ and $\kappa=5$.}
\label{bgb_fig1}                      
\end{center}
\end{figure}

\begin{figure}[H]
\begin{center}
\raisebox{0ex}{\includegraphics[width=0.35\linewidth]{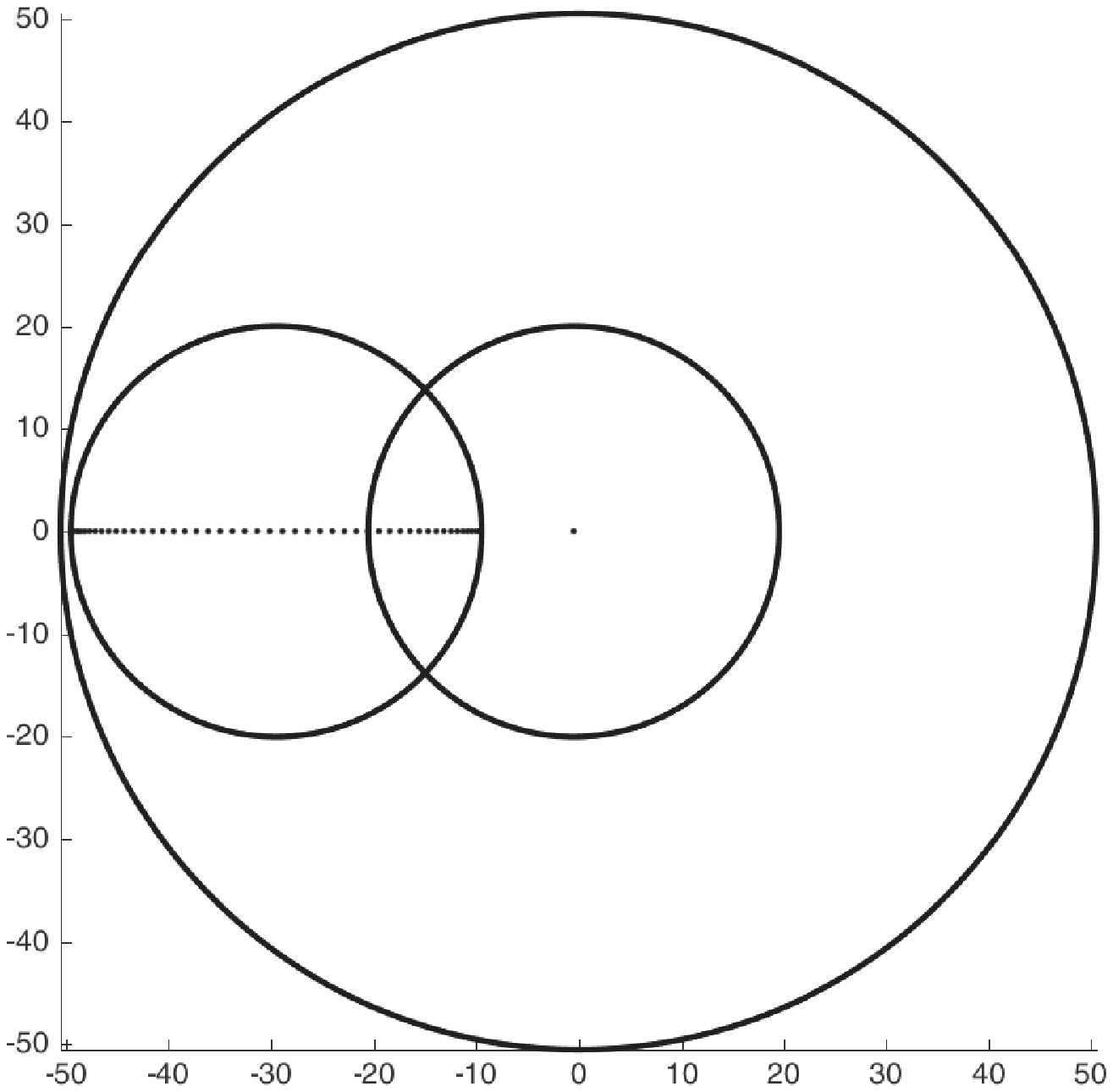}}
\hskip 2cm
\raisebox{0ex}{\includegraphics[width=0.35\linewidth]{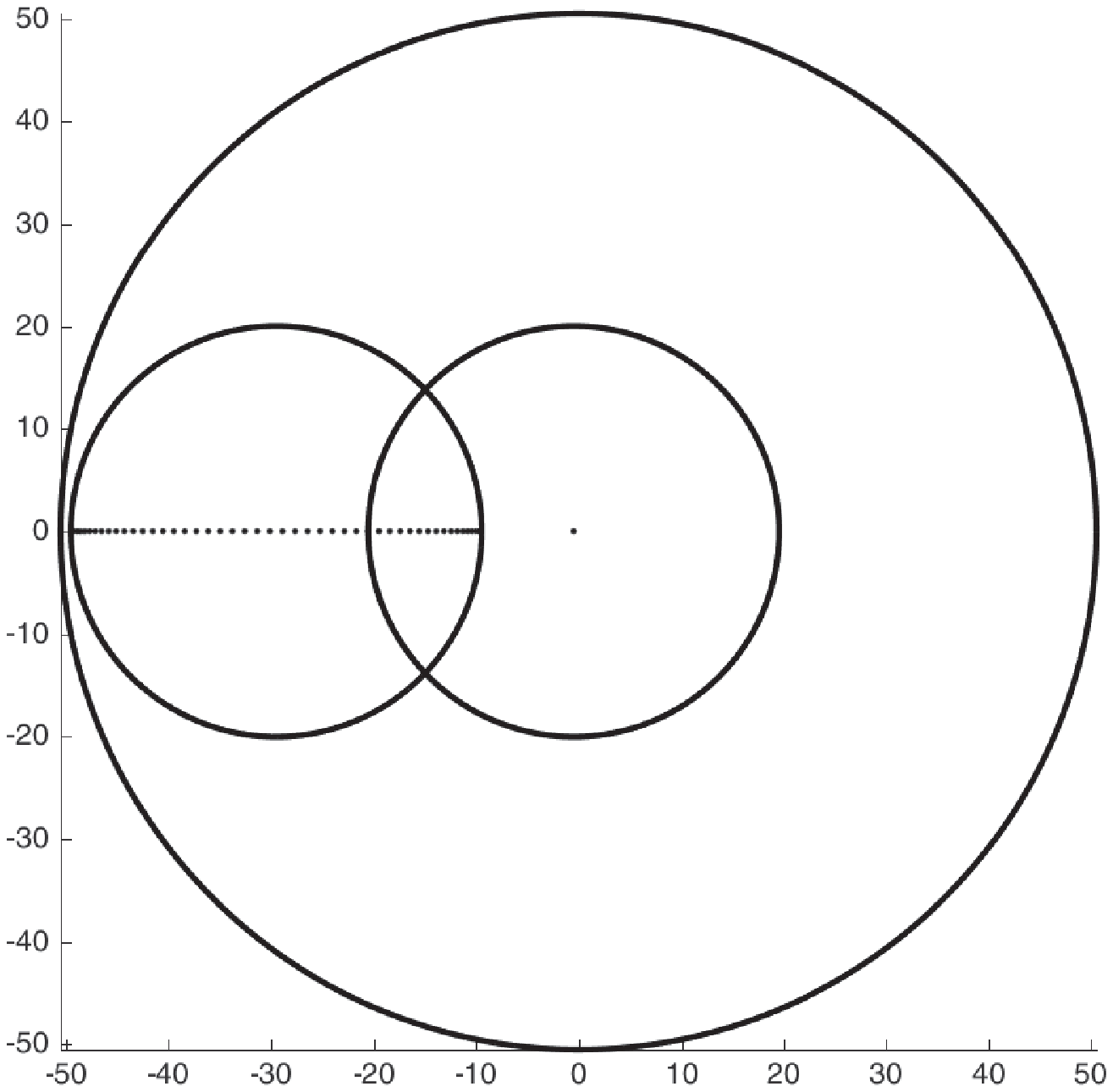}}
\caption{Inclusion regions for Example~1 with $\tau=10$ and $\kappa=5$.}
\label{bgb_fig2}                      
\end{center}
\end{figure}

\begin{figure}[H]
\begin{center}
\raisebox{0ex}{\includegraphics[width=0.35\linewidth]{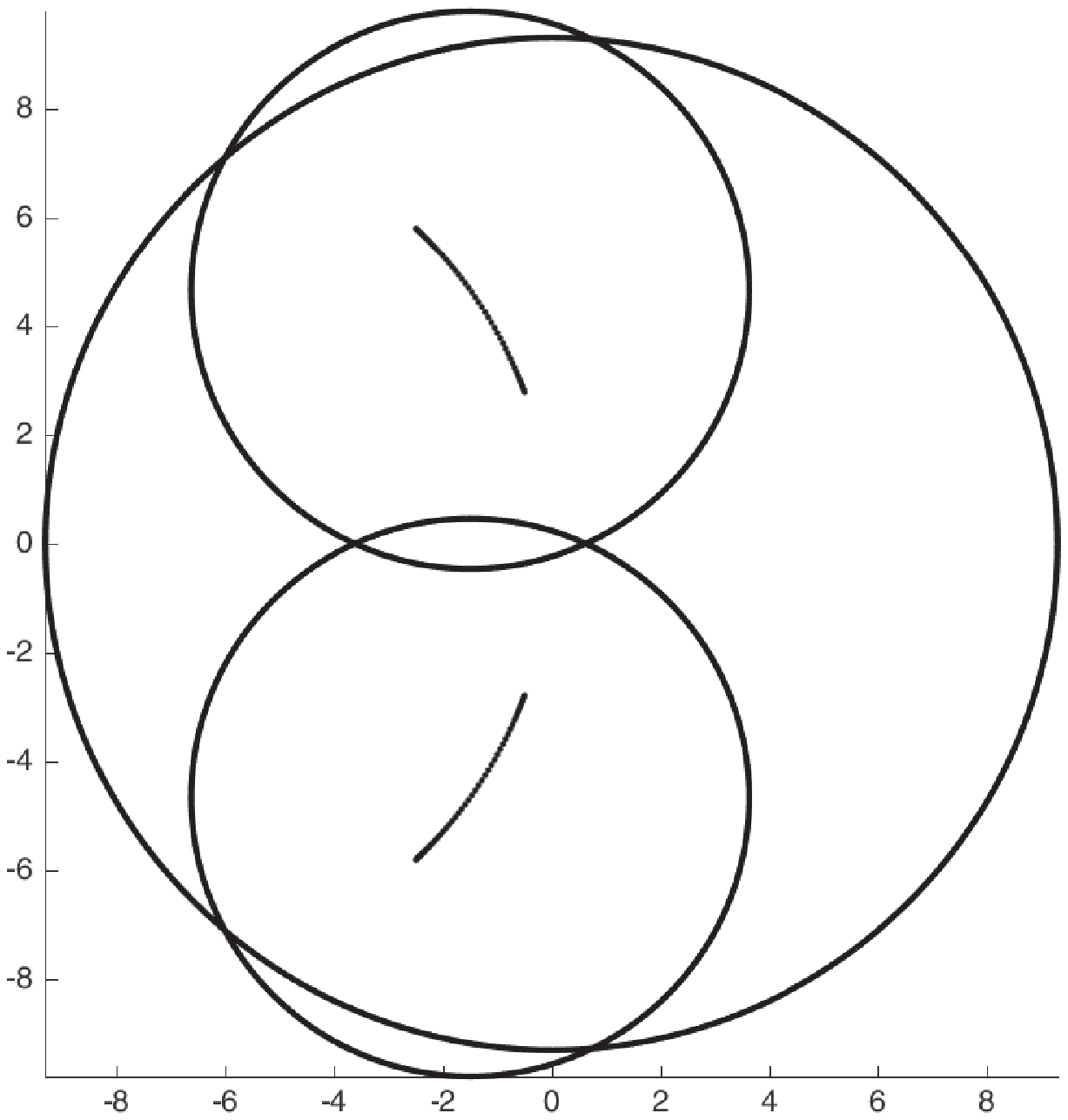}}
\hskip 2cm
\raisebox{0ex}{\includegraphics[width=0.35\linewidth]{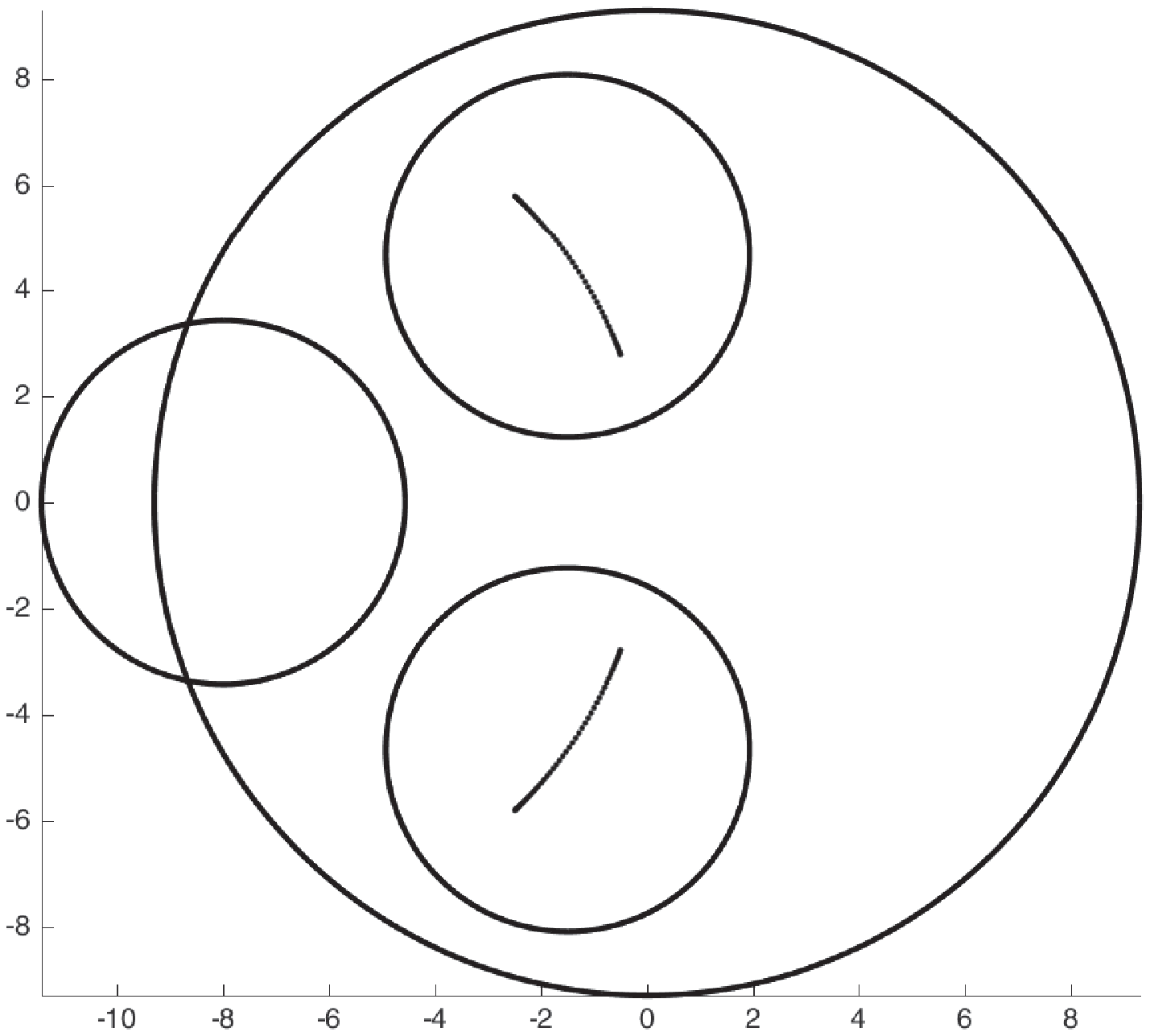}}
\caption{Inclusion regions for Example~1 with $\tau=1$ and $\kappa=8$.}
\label{bgb_fig3}                      
\end{center}
\end{figure}

\begin{figure}[H]
\begin{center}
\raisebox{0ex}{\includegraphics[width=0.35\linewidth]{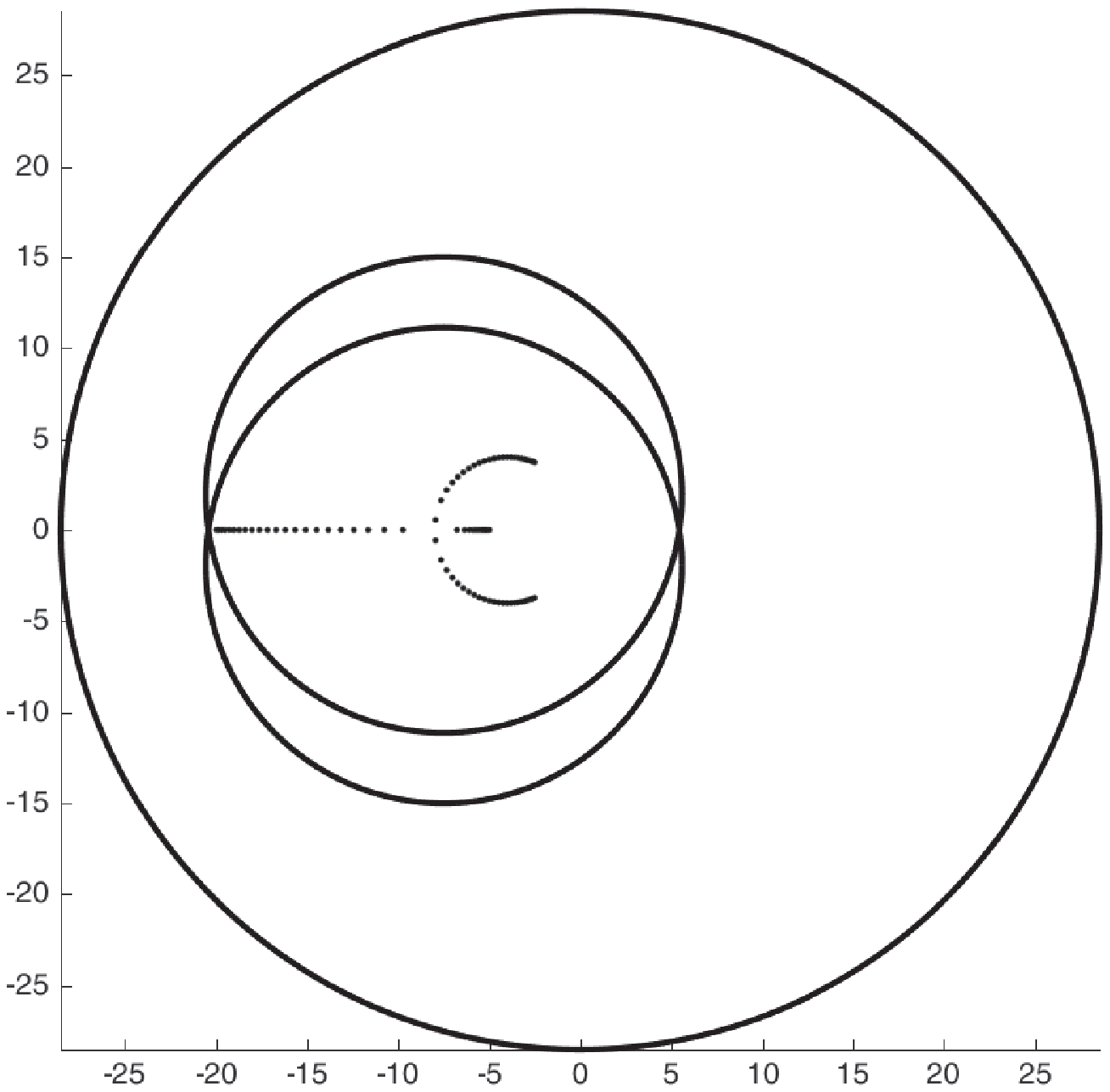}}
\hskip 2cm
\raisebox{0ex}{\includegraphics[width=0.35\linewidth]{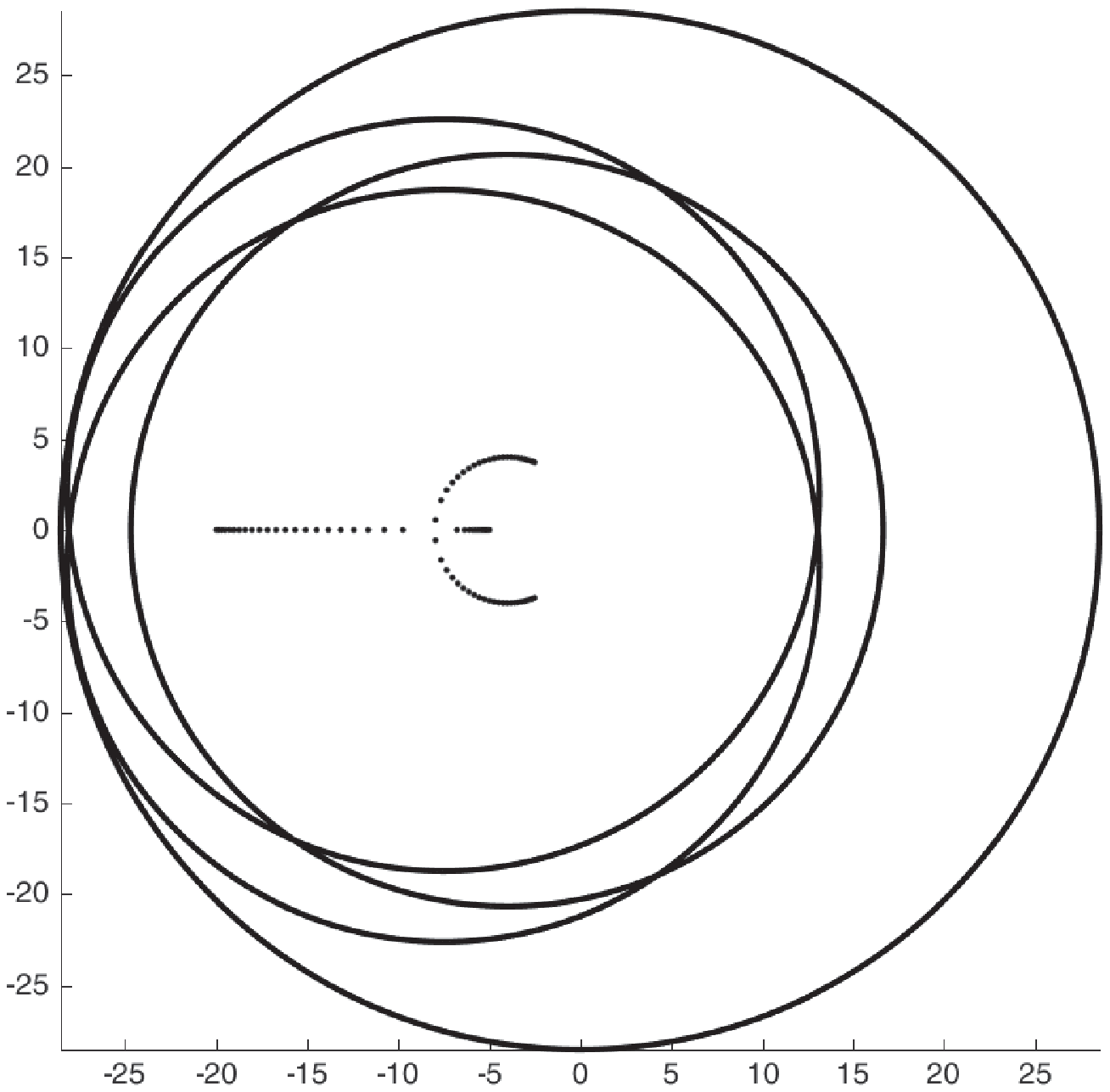}}
\caption{Inclusion regions for Example~1 with $\tau=5$ and $\kappa=20$.}
\label{bgb_fig4}                      
\end{center}
\end{figure}

\begin{figure}[H]
\begin{center}
\raisebox{0ex}{\includegraphics[width=0.35\linewidth]{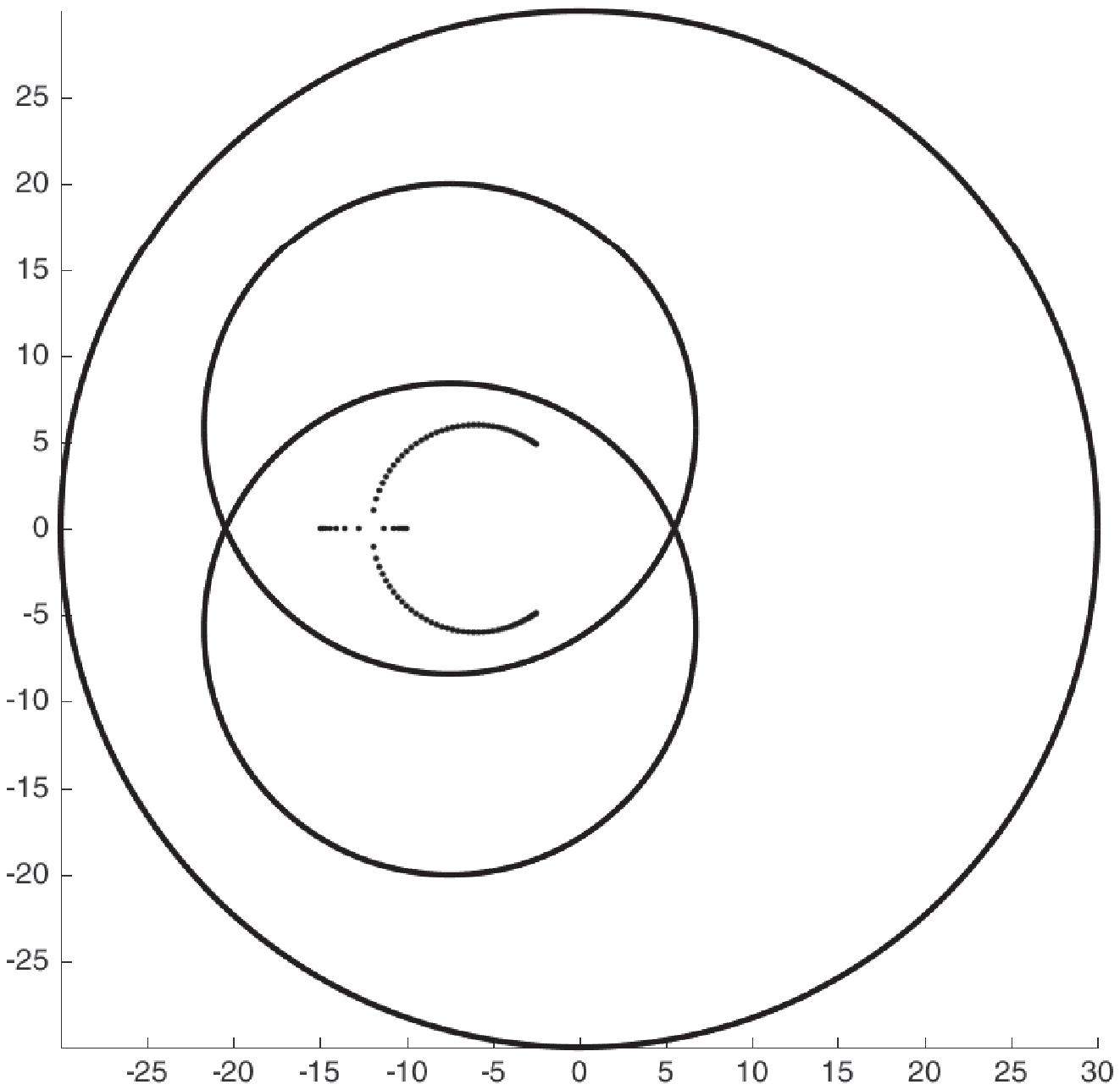}}
\hskip 2cm
\raisebox{0ex}{\includegraphics[width=0.35\linewidth]{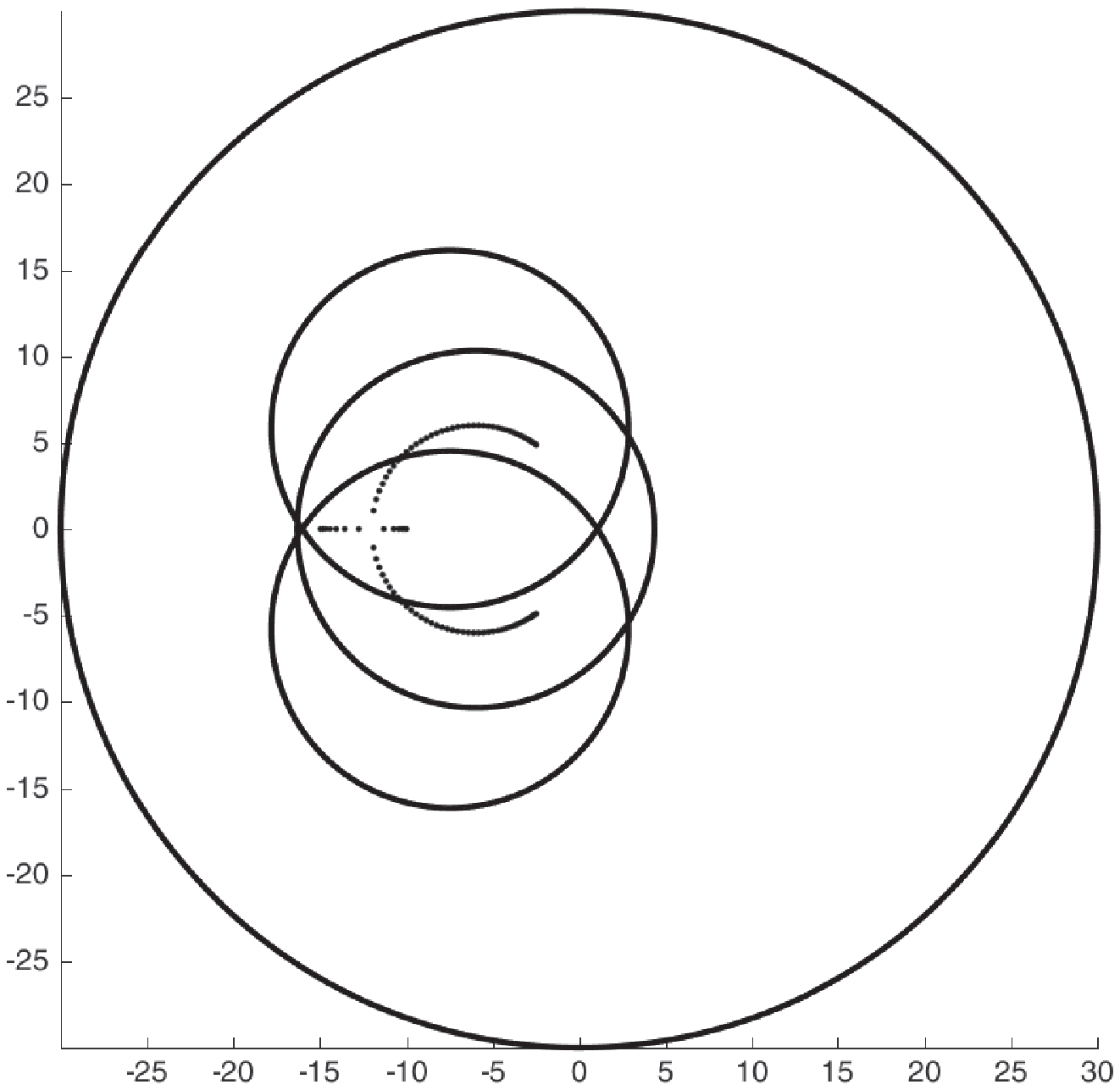}}
\caption{Inclusion regions for Example~1 with $\tau=5$ and $\kappa=30$.}
\label{bgb_fig5}                      
\end{center}
\end{figure}

The following figures show the eigenvalue inclusion regions for a few representative
values of $\tau$ and $\kappa$. All eigenvalues have negative real parts since the coefficient matrices are all strictly positive definite
(see \cite{TisseurMeerbergen}).
In each figure, the large circle centered at the origin is the circle obtained from Theorem~\ref{genCauchy}, namely, Cauchy's theorem for matrix polynomials;
its radius is the Cauchy radius of $P$ and we will refer to it as the \emph{Cauchy disk} of $P$.
On the left, the smaller disks represent the inclusion region obtained from Theorem~\ref{maintheorem} for the Newton
basis, while those on the right are for the basis $\mathcal{B}$.
Figure~\ref{bgb_fig1} and Figure~\ref{bgb_fig2} show the eigenvalue inclusion regions for $\tau=3,\kappa=5$ and $\tau=10,\kappa=5$, respectively,
which are the values used in~\cite{TisseurMeerbergen}.
The dots are the eigenvalues, which are added for reference. 
For Figure~\ref{bgb_fig2}, \ref{bgb_fig3}, \ref{bgb_fig4}, and ~\ref{bgb_fig5}, 
the values for the pair $(\tau,\kappa)$ are $(1,8)$, $(5,20)$, $(5,30)$, and $(5,80)$, respectively.
When $\tau$ is large relative to $\kappa$, the inclusion regions are almost identical for both bases, as for the $(10,5)$ case.

These figures clearly show that using a more general basis can significantly reduce the eigenvalue inclusion regions, when compared to the disk obtained
from Theorem~\ref{genCauchy}, which is often the best one can obtain for the power basis. Sometimes the Newton basis is better, and sometimes it is the 
more general basis $\mathcal{B}$ that produces the smaller inclusion region. Of special interest is Figure~\ref{bgb_fig3}, where the eigenvalues are split
among the top and bottom disks, $50$ in each disk, as predicted by the theorem, since the middle disk does not contain any zeros of $\qtwo$.
We remark that it would not be possible to obtain such an inclusion region from any of the classical bounds.
%
% Figures created with bgb_example1.m 
%
%%%%%%%%%%%%%%%%%%%%%%%%%%%%%%%%%%%%%%%%%%%%%%%%%%%%%%%%%%%%%%%%%%%%%%%%%%%%%%%%%%%%%%%%%%%%%%%%%%%%%%%%%%%%%%%%%%%%%%%%%%%%%%%%%%%%%%%%%%%%%%%%%%%%%%%%%%%%
%%%%%%%%%%%%%%%%%%%%%%%%%%%%%%%%%%%%%%%%%%%%%%%%%%%%%%%%%%%%%%%%%%%%%%%%%%%%%%%%%%%%%%%%%%%%%%%%%%%%%%%%%%%%%%%%%%%%%%%%%%%%%%%%%%%%%%%%%%%%%%%%%%%%%%%%%%%%
%
%                                           F I G U R E  1
%
%%%%%%%%%%%%%%%%%%%%%%%%%%%%%%%%%%%%%%%%%%%%%%%%%%%%%%%%%%%%%%%%%%%%%%%%%%%%%%%%%%%%%%%%%%%%%%%%%%%%%%%%%%%%%%%%%%%%%%%%%%%%%%%%%%%%%%%%%%%%%%%%%%%%%%%%%%%%
%%%%%%%%%%%%%%%%%%%%%%%%%%%%%%%%%%%%%%%%%%%%%%%%%%%%%%%%%%%%%%%%%%%%%%%%%%%%%%%%%%%%%%%%%%%%%%%%%%%%%%%%%%%%%%%%%%%%%%%%%%%%%%%%%%%%%%%%%%%%%%%%%%%%%%%%%%%%
%%%%%%%%%%%%%%%%%%%%%%%%%%%%%%%%%%%%%%%%%%%%%%%%%%%%%%%%%%%%%%%%%%%%%%%%%%%%%%%%%%%%%%%%%%%%%%%%%%%%%%%%%%%%%%%%%%%%%%%%%%%%%%%%%%%%%%%%%%%%%%%%%%%%%%%%%%%%
%%%%%%%%%%%%%%%%%%%%%%%%%%%%%%%%%%%%%%%%%%%%%%%%%%%%%%%%%%%%%%%%%%%%%%%%%%%%%%%%%%%%%%%%%%%%%%%%%%%%%%%%%%%%%%%%%%%%%%%%%%%%%%%%%%%%%%%%%%%%%%%%%%%%%%%%%%%%

%%%%%%%%%%%%%%%%%%%%%%%%%%%%%%%%%%%%%%%%%%%%%%%%%%%%%%%%%%%%%%%%%%%%%%%%%%%%%%%%%%%%%%%%%%%%%%%%%%%%%%%%%%%%%%%%%%%%%%%%%%%%%%%%%%%%%%%%%%%%%%%%%%%%%%%%%%%%
%%%%%%%%%%%%%%%%%%%%%%%%%%%%%%%%%%%%%%%%%%%%%%%%%%%%%%%%%%%%%%%%%%%%%%%%%%%%%%%%%%%%%%%%%%%%%%%%%%%%%%%%%%%%%%%%%%%%%%%%%%%%%%%%%%%%%%%%%%%%%%%%%%%%%%%%%%%%
%
%                                           F I G U R E  2
%
%%%%%%%%%%%%%%%%%%%%%%%%%%%%%%%%%%%%%%%%%%%%%%%%%%%%%%%%%%%%%%%%%%%%%%%%%%%%%%%%%%%%%%%%%%%%%%%%%%%%%%%%%%%%%%%%%%%%%%%%%%%%%%%%%%%%%%%%%%%%%%%%%%%%%%%%%%%%
%%%%%%%%%%%%%%%%%%%%%%%%%%%%%%%%%%%%%%%%%%%%%%%%%%%%%%%%%%%%%%%%%%%%%%%%%%%%%%%%%%%%%%%%%%%%%%%%%%%%%%%%%%%%%%%%%%%%%%%%%%%%%%%%%%%%%%%%%%%%%%%%%%%%%%%%%%%%
%%%%%%%%%%%%%%%%%%%%%%%%%%%%%%%%%%%%%%%%%%%%%%%%%%%%%%%%%%%%%%%%%%%%%%%%%%%%%%%%%%%%%%%%%%%%%%%%%%%%%%%%%%%%%%%%%%%%%%%%%%%%%%%%%%%%%%%%%%%%%%%%%%%%%%%%%%%%
%%%%%%%%%%%%%%%%%%%%%%%%%%%%%%%%%%%%%%%%%%%%%%%%%%%%%%%%%%%%%%%%%%%%%%%%%%%%%%%%%%%%%%%%%%%%%%%%%%%%%%%%%%%%%%%%%%%%%%%%%%%%%%%%%%%%%%%%%%%%%%%%%%%%%%%%%%%%

%%%%%%%%%%%%%%%%%%%%%%%%%%%%%%%%%%%%%%%%%%%%%%%%%%%%%%%%%%%%%%%%%%%%%%%%%%%%%%%%%%%%%%%%%%%%%%%%%%%%%%%%%%%%%%%%%%%%%%%%%%%%%%%%%%%%%%%%%%%%%%%%%%%%%%%%%%%%
%%%%%%%%%%%%%%%%%%%%%%%%%%%%%%%%%%%%%%%%%%%%%%%%%%%%%%%%%%%%%%%%%%%%%%%%%%%%%%%%%%%%%%%%%%%%%%%%%%%%%%%%%%%%%%%%%%%%%%%%%%%%%%%%%%%%%%%%%%%%%%%%%%%%%%%%%%%%
%
%                                           F I G U R E  3
%
%%%%%%%%%%%%%%%%%%%%%%%%%%%%%%%%%%%%%%%%%%%%%%%%%%%%%%%%%%%%%%%%%%%%%%%%%%%%%%%%%%%%%%%%%%%%%%%%%%%%%%%%%%%%%%%%%%%%%%%%%%%%%%%%%%%%%%%%%%%%%%%%%%%%%%%%%%%%
%%%%%%%%%%%%%%%%%%%%%%%%%%%%%%%%%%%%%%%%%%%%%%%%%%%%%%%%%%%%%%%%%%%%%%%%%%%%%%%%%%%%%%%%%%%%%%%%%%%%%%%%%%%%%%%%%%%%%%%%%%%%%%%%%%%%%%%%%%%%%%%%%%%%%%%%%%%%
%%%%%%%%%%%%%%%%%%%%%%%%%%%%%%%%%%%%%%%%%%%%%%%%%%%%%%%%%%%%%%%%%%%%%%%%%%%%%%%%%%%%%%%%%%%%%%%%%%%%%%%%%%%%%%%%%%%%%%%%%%%%%%%%%%%%%%%%%%%%%%%%%%%%%%%%%%%%
%%%%%%%%%%%%%%%%%%%%%%%%%%%%%%%%%%%%%%%%%%%%%%%%%%%%%%%%%%%%%%%%%%%%%%%%%%%%%%%%%%%%%%%%%%%%%%%%%%%%%%%%%%%%%%%%%%%%%%%%%%%%%%%%%%%%%%%%%%%%%%%%%%%%%%%%%%%%

%%%%%%%%%%%%%%%%%%%%%%%%%%%%%%%%%%%%%%%%%%%%%%%%%%%%%%%%%%%%%%%%%%%%%%%%%%%%%%%%%%%%%%%%%%%%%%%%%%%%%%%%%%%%%%%%%%%%%%%%%%%%%%%%%%%%%%%%%%%%%%%%%%%%%%%%%%%%
%%%%%%%%%%%%%%%%%%%%%%%%%%%%%%%%%%%%%%%%%%%%%%%%%%%%%%%%%%%%%%%%%%%%%%%%%%%%%%%%%%%%%%%%%%%%%%%%%%%%%%%%%%%%%%%%%%%%%%%%%%%%%%%%%%%%%%%%%%%%%%%%%%%%%%%%%%%%
%
%                                           F I G U R E  4
%
%%%%%%%%%%%%%%%%%%%%%%%%%%%%%%%%%%%%%%%%%%%%%%%%%%%%%%%%%%%%%%%%%%%%%%%%%%%%%%%%%%%%%%%%%%%%%%%%%%%%%%%%%%%%%%%%%%%%%%%%%%%%%%%%%%%%%%%%%%%%%%%%%%%%%%%%%%%%
%%%%%%%%%%%%%%%%%%%%%%%%%%%%%%%%%%%%%%%%%%%%%%%%%%%%%%%%%%%%%%%%%%%%%%%%%%%%%%%%%%%%%%%%%%%%%%%%%%%%%%%%%%%%%%%%%%%%%%%%%%%%%%%%%%%%%%%%%%%%%%%%%%%%%%%%%%%%
%%%%%%%%%%%%%%%%%%%%%%%%%%%%%%%%%%%%%%%%%%%%%%%%%%%%%%%%%%%%%%%%%%%%%%%%%%%%%%%%%%%%%%%%%%%%%%%%%%%%%%%%%%%%%%%%%%%%%%%%%%%%%%%%%%%%%%%%%%%%%%%%%%%%%%%%%%%%
%%%%%%%%%%%%%%%%%%%%%%%%%%%%%%%%%%%%%%%%%%%%%%%%%%%%%%%%%%%%%%%%%%%%%%%%%%%%%%%%%%%%%%%%%%%%%%%%%%%%%%%%%%%%%%%%%%%%%%%%%%%%%%%%%%%%%%%%%%%%%%%%%%%%%%%%%%%%

%%%%%%%%%%%%%%%%%%%%%%%%%%%%%%%%%%%%%%%%%%%%%%%%%%%%%%%%%%%%%%%%%%%%%%%%%%%%%%%%%%%%%%%%%%%%%%%%%%%%%%%%%%%%%%%%%%%%%%%%%%%%%%%%%%%%%%%%%%%%%%%%%%%%%%%%%%%%
%%%%%%%%%%%%%%%%%%%%%%%%%%%%%%%%%%%%%%%%%%%%%%%%%%%%%%%%%%%%%%%%%%%%%%%%%%%%%%%%%%%%%%%%%%%%%%%%%%%%%%%%%%%%%%%%%%%%%%%%%%%%%%%%%%%%%%%%%%%%%%%%%%%%%%%%%%%%
%
%                                           F I G U R E  5
%
%%%%%%%%%%%%%%%%%%%%%%%%%%%%%%%%%%%%%%%%%%%%%%%%%%%%%%%%%%%%%%%%%%%%%%%%%%%%%%%%%%%%%%%%%%%%%%%%%%%%%%%%%%%%%%%%%%%%%%%%%%%%%%%%%%%%%%%%%%%%%%%%%%%%%%%%%%%%
%%%%%%%%%%%%%%%%%%%%%%%%%%%%%%%%%%%%%%%%%%%%%%%%%%%%%%%%%%%%%%%%%%%%%%%%%%%%%%%%%%%%%%%%%%%%%%%%%%%%%%%%%%%%%%%%%%%%%%%%%%%%%%%%%%%%%%%%%%%%%%%%%%%%%%%%%%%%
%%%%%%%%%%%%%%%%%%%%%%%%%%%%%%%%%%%%%%%%%%%%%%%%%%%%%%%%%%%%%%%%%%%%%%%%%%%%%%%%%%%%%%%%%%%%%%%%%%%%%%%%%%%%%%%%%%%%%%%%%%%%%%%%%%%%%%%%%%%%%%%%%%%%%%%%%%%%
%%%%%%%%%%%%%%%%%%%%%%%%%%%%%%%%%%%%%%%%%%%%%%%%%%%%%%%%%%%%%%%%%%%%%%%%%%%%%%%%%%%%%%%%%%%%%%%%%%%%%%%%%%%%%%%%%%%%%%%%%%%%%%%%%%%%%%%%%%%%%%%%%%%%%%%%%%%%

%%%%%%%%%%%%%%%%%%%%%%%%%%%%%%%%%%%%%%%%%%%%%%%%%%%%%%%%%%%%%%%%%%%%%%%%%%%%%%%%%%%%%%%%%%%%%%%%%%%%%%%%%%%%%%%%%%%%%%%%%%%%%%%%%%%%%%%%%%%%%%%%%%%%%%%%%%%%
%%%%%%%%%%%%%%%%%%%%%%%%%%%%%%%%%%%%%%%%%%%%%%%%%%%%%%%%%%%%%%%%%%%%%%%%%%%%%%%%%%%%%%%%%%%%%%%%%%%%%%%%%%%%%%%%%%%%%%%%%%%%%%%%%%%%%%%%%%%%%%%%%%%%%%%%%%%%
%
%                                           F I G U R E  6
%
%%%%%%%%%%%%%%%%%%%%%%%%%%%%%%%%%%%%%%%%%%%%%%%%%%%%%%%%%%%%%%%%%%%%%%%%%%%%%%%%%%%%%%%%%%%%%%%%%%%%%%%%%%%%%%%%%%%%%%%%%%%%%%%%%%%%%%%%%%%%%%%%%%%%%%%%%%%%
%%%%%%%%%%%%%%%%%%%%%%%%%%%%%%%%%%%%%%%%%%%%%%%%%%%%%%%%%%%%%%%%%%%%%%%%%%%%%%%%%%%%%%%%%%%%%%%%%%%%%%%%%%%%%%%%%%%%%%%%%%%%%%%%%%%%%%%%%%%%%%%%%%%%%%%%%%%%
\begin{figure}[H]
\begin{center}
\raisebox{0ex}{\includegraphics[width=0.30\linewidth]{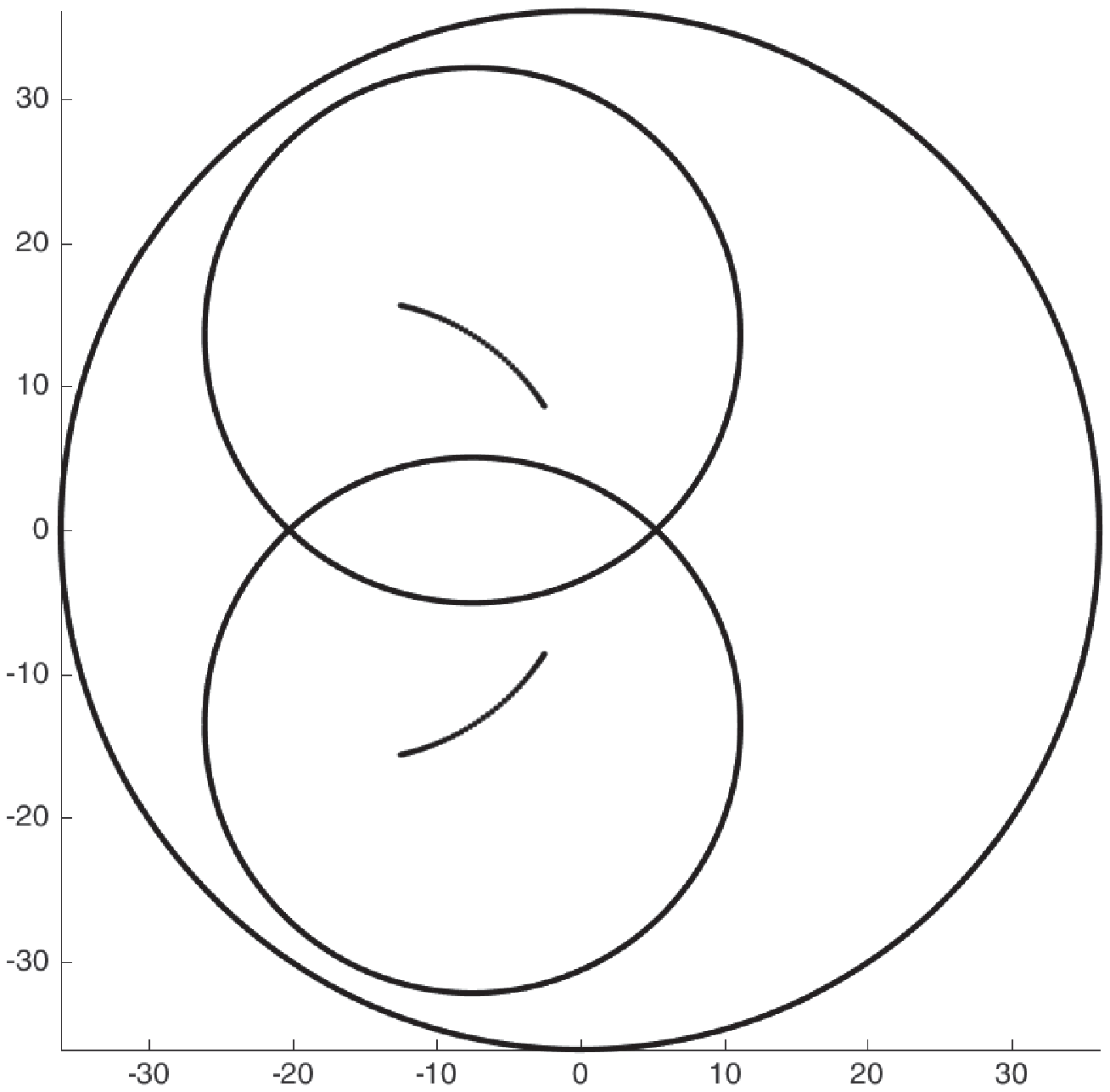}}
\hskip 2cm
\raisebox{0ex}{\includegraphics[width=0.30\linewidth]{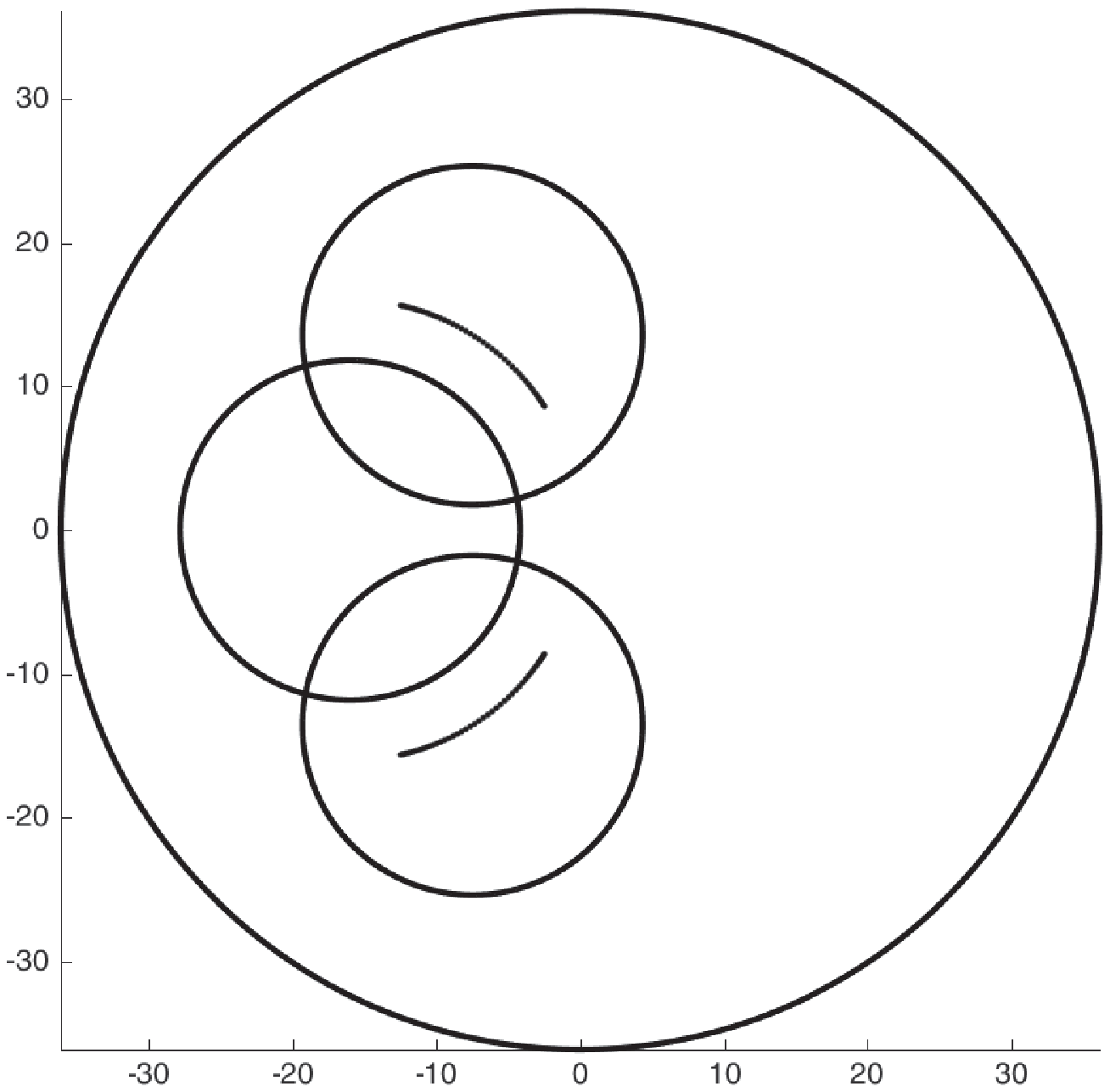}}
\caption{Inclusion regions for Example~1 with $\tau=5$ and $\kappa=80$.}
\label{bgb_fig6}                      
\end{center}
\end{figure}
%%%%%%%%%%%%%%%%%%%%%%%%%%%%%%%%%%%%%%%%%%%%%%%%%%%%%%%%%%%%%%%%%%%%%%%%%%%%%%%%%%%%%%%%%%%%%%%%%%%%%%%%%%%%%%%%%%%%%%%%%%%%%%%%%%%%%%%%%%%%%%%%%%%%%%%%%%%%
%%%%%%%%%%%%%%%%%%%%%%%%%%%%%%%%%%%%%%%%%%%%%%%%%%%%%%%%%%%%%%%%%%%%%%%%%%%%%%%%%%%%%%%%%%%%%%%%%%%%%%%%%%%%%%%%%%%%%%%%%%%%%%%%%%%%%%%%%%%%%%%%%%%%%%%%%%%%

%%%%%%%%%%%%%%%%%%%%%%%%%%%%%%%%%%%%%%%%%%%%%%%%%%%%%%%%%%%%%%%%%%%%%%%%%%%%%%%%%%%%%%%%%%%%%%%%%%%%%%%%%%%%%%%%%%%%%%%%%%%%%%%%%%%%%%%%%%%%%%%%%%%%%%%%%%%%%%%%%%

%%%%%%%%%%%%%%%%%%%%%%%%%%%%%%%%%%%%%%%%%%%%%%%%%%%%%%%%%%%%%%%%%%%%%%%%%%%%%%%%%%%%%%%%%%%%%%%%%%%%%%%%%%%%%%%%%%%%%%%%%%%%%%%%%%%%%%%%%%%%%%%%%%%%%%%%%%%%%%%%%%
%
%
%                                                         EXAMPLE 2                                                               
%
%
%%%%%%%%%%%%%%%%%%%%%%%%%%%%%%%%%%%%%%%%%%%%%%%%%%%%%%%%%%%%%%%%%%%%%%%%%%%%%%%%%%%%%%%%%%%%%%%%%%%%%%%%%%%%%%%%%%%%%%%%%%%%%%%%%%%%%%%%%%%%%%%%%%%%%%%%%%%%%%%%%%
\noindent {\bf Example 2.} In this example from~\cite{BHMST} and~\cite{CCVG}, we consider a quadratic polynomial produced by the finite-element discretization of a 
time-harmonic wave equation for the acoustic pressure on the unit square $[0,1] \times [0,1]$. The eigenvalues lie in the upper half of the complex plane.
Here we have $m=\ell(\ell-1)$, where $\ell=1/h$ and $h$ is the mesh size. Defining the $\ell \times \ell$
matrix $S_{\ell}$ and the $(\ell-1)\times (\ell-1)$ matrix $T_{\ell-1}$ as 
\bdis
S_{\ell} = 
\begin{pmatrix}
4    &    -1     &           &      \\
-1   &  \ddots   &  \ddots   &      \\
     &  \ddots   &    4      &  -1  \\
     &           &   -1      &   2  \\
\end{pmatrix}
\;\; , \;\;
\text{and} \;\; 
T_{\ell-1} = 
\begin{pmatrix}
0    &    -1     &           &      \\
-1   &  \ddots   &  \ddots   &      \\
     &  \ddots   &  \ddots   &  -1  \\
     &           &   -1      &   0  \\
\end{pmatrix}
\;\; , 
\edis
the coefficients of the quadratic matrix polynomial $P(z)=\Atwo z^{2} + \Aone z + \Azero$ are given by
\bdis
A_{0} = I_{\ell-1} \otms \, S_{\ell} + T_{\ell-1} \otms \, \lb -I_{\ell} + \dfrac{1}{2} e_{\ell}e_{\ell}^{T} \rb 
\;\; , \;\;
A_{1} = \dfrac{2\pi}{\ell \zeta} I_{\ell-1} \otms \, e_{\ell}e_{\ell}^{T} 
\;\; , \;\;
\edis
\bdis
A_{2} = -\dfrac{4\pi^{2}}{\ell^{2}} I_{\ell-1} \otms \lb I_{\ell} - \dfrac{1}{2} e_{\ell}e_{\ell}^{T} \rb 
\; ,
\edis
where the complex number $\zeta$ is the impedance, $e_{i}$ is the $i$th standard unit vector,
and the Kronecker product of two matrices $A \otms B$ is the block matrix $(a_{ij}B)$.
Since $A_{2}$ is nonsingular and diagonal, it is an easy matter to compute $\Atwo^{-1}P=Iz^{2}+B_{1}z+B_{0}$, where $B_{1}=\Atwo^{-1}\Aone$ 
and $B_{0}=\Atwo^{-1}\Azero$. The matrix $B_{0}$ is diagonally dominant for most of its rows and columns.
We can now conveniently use the results from Subsection~\ref{quadraticproblems} to express $P$ in the basis $\mathcal{N}$.
In this example, as in the next, we will only consider the Newton basis, since the disks in both bases are not significantly different in size.

The diagonals are not constant, and to minimize their $1$-norm we use the observation about the minimization problem in~(\ref{minproblem}).
We will once again aim to choose nodes that minimize the $1$-norms of 
the coefficients of $P$. We set $\text{diag} \bigl ( B_{j} \bigr ) = C_{j} + i D_{j}$ for $j=0,1$, and, 
in light of the above observation about~(\ref{minproblem}), we define
\bdis
\mu_{j} = \dfrac{1}{2} \bigl ( \min \lb C_{j} \rb + \max \lb C_{j} \rb \bigr ) + \dfrac{i}{2} \bigl ( \min \lb D_{j} \rb + \max \lb D_{j} \rb  \bigr )  
\qquad \text{($j=0,1$)} \; . 
\edis
From the expression in~(\ref{pinbases}) for the matrix coefficients of $\fone$, we see that a reasonable choice for the nodes is to choose them so that
$a+b = -\mu_{1}$. Since from~(\ref{pinbases}) the coefficient of $\fzero$ can be written as 
\bdis
aB_{1} + B_{0} + a^{2} I = a^{2}I +  \mu_{1} a I + \mu_{0}I +  a \bigl ( B_{1} - \mu_{1}I \bigr ) + \bigl ( B_{0}-\mu_{0}I \bigr ) \; ,
\edis
% where $B_{0} - \mu_{0}$ has zero diagonal, 
we choose $a$ as that solution of $a^{2}+\mu_{1} a +\mu_{0} = 0$ that minimizes $\|aB_{1} + B_{0} + a^{2} I\|_{1}$.

\begin{figure}[H]
\begin{center}
\raisebox{0ex}{\includegraphics[width=0.45\linewidth]{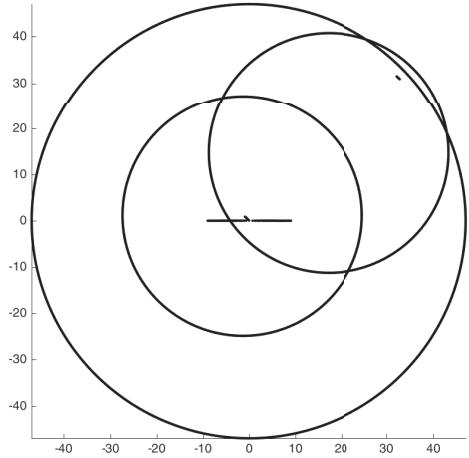}}
\caption{Inclusion regions for Example~2 with $\zeta=0.1 + 0.1i$.}
\label{bgb_fig7}                      
\end{center}
\end{figure}

\begin{figure}[H]
\begin{center}
\raisebox{0ex}{\includegraphics[width=0.45\linewidth]{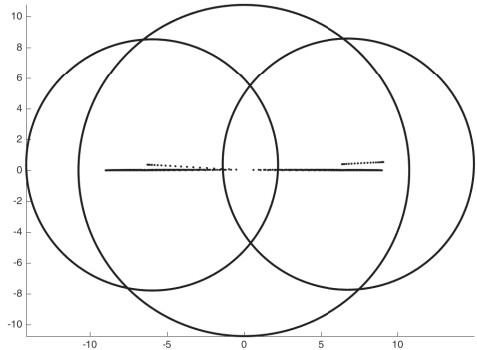}}
\caption{Inclusion regions for Example~2 with $\zeta=2+2i$.}
\label{bgb_fig8}                      
\end{center}
\end{figure}

Figures~\ref{bgb_fig7} and~\ref{bgb_fig8} show the eigenvalue inclusion regions for $\zeta= 0.1 + 0.1i$ and $2+2i$, respectively, with $h=0.05$, so
that the matrix coefficients are of size $380 \times 380$. 
The large circle centered at the origin is, as before, the Cauchy disk of $P$, while the black dots are the eigenvalues. 

For very small or very large values of $|\zeta |$, the disks in the Newton basis are not significantly different from the ones in Figure~\ref{bgb_fig7} 
and Figure~\ref{bgb_fig8}, respectively.
From these results, it is clear that using a generalized basis here clearly allows a significant part of the Cauchy disk to be discarded as a possible location 
for the eigenvalues.
%
% Figures created with bgb_example2.m 
%
%%%%%%%%%%%%%%%%%%%%%%%%%%%%%%%%%%%%%%%%%%%%%%%%%%%%%%%%%%%%%%%%%%%%%%%%%%%%%%%%%%%%%%%%%%%%%%%%%%%%%%%%%%%%%%%%%%%%%%%%%%%%%%%%%%%%%%%%%%%%%%%%%%%%%%%%%%%%
%%%%%%%%%%%%%%%%%%%%%%%%%%%%%%%%%%%%%%%%%%%%%%%%%%%%%%%%%%%%%%%%%%%%%%%%%%%%%%%%%%%%%%%%%%%%%%%%%%%%%%%%%%%%%%%%%%%%%%%%%%%%%%%%%%%%%%%%%%%%%%%%%%%%%%%%%%%%
%
%                                           F I G U R E  7
%
%%%%%%%%%%%%%%%%%%%%%%%%%%%%%%%%%%%%%%%%%%%%%%%%%%%%%%%%%%%%%%%%%%%%%%%%%%%%%%%%%%%%%%%%%%%%%%%%%%%%%%%%%%%%%%%%%%%%%%%%%%%%%%%%%%%%%%%%%%%%%%%%%%%%%%%%%%%%
%%%%%%%%%%%%%%%%%%%%%%%%%%%%%%%%%%%%%%%%%%%%%%%%%%%%%%%%%%%%%%%%%%%%%%%%%%%%%%%%%%%%%%%%%%%%%%%%%%%%%%%%%%%%%%%%%%%%%%%%%%%%%%%%%%%%%%%%%%%%%%%%%%%%%%%%%%%%
%%%%%%%%%%%%%%%%%%%%%%%%%%%%%%%%%%%%%%%%%%%%%%%%%%%%%%%%%%%%%%%%%%%%%%%%%%%%%%%%%%%%%%%%%%%%%%%%%%%%%%%%%%%%%%%%%%%%%%%%%%%%%%%%%%%%%%%%%%%%%%%%%%%%%%%%%%%%
%%%%%%%%%%%%%%%%%%%%%%%%%%%%%%%%%%%%%%%%%%%%%%%%%%%%%%%%%%%%%%%%%%%%%%%%%%%%%%%%%%%%%%%%%%%%%%%%%%%%%%%%%%%%%%%%%%%%%%%%%%%%%%%%%%%%%%%%%%%%%%%%%%%%%%%%%%%%

%%%%%%%%%%%%%%%%%%%%%%%%%%%%%%%%%%%%%%%%%%%%%%%%%%%%%%%%%%%%%%%%%%%%%%%%%%%%%%%%%%%%%%%%%%%%%%%%%%%%%%%%%%%%%%%%%%%%%%%%%%%%%%%%%%%%%%%%%%%%%%%%%%%%%%%%%%%%
%%%%%%%%%%%%%%%%%%%%%%%%%%%%%%%%%%%%%%%%%%%%%%%%%%%%%%%%%%%%%%%%%%%%%%%%%%%%%%%%%%%%%%%%%%%%%%%%%%%%%%%%%%%%%%%%%%%%%%%%%%%%%%%%%%%%%%%%%%%%%%%%%%%%%%%%%%%%
%
%                                           F I G U R E  8
%
%%%%%%%%%%%%%%%%%%%%%%%%%%%%%%%%%%%%%%%%%%%%%%%%%%%%%%%%%%%%%%%%%%%%%%%%%%%%%%%%%%%%%%%%%%%%%%%%%%%%%%%%%%%%%%%%%%%%%%%%%%%%%%%%%%%%%%%%%%%%%%%%%%%%%%%%%%%%
%%%%%%%%%%%%%%%%%%%%%%%%%%%%%%%%%%%%%%%%%%%%%%%%%%%%%%%%%%%%%%%%%%%%%%%%%%%%%%%%%%%%%%%%%%%%%%%%%%%%%%%%%%%%%%%%%%%%%%%%%%%%%%%%%%%%%%%%%%%%%%%%%%%%%%%%%%%%
%%%%%%%%%%%%%%%%%%%%%%%%%%%%%%%%%%%%%%%%%%%%%%%%%%%%%%%%%%%%%%%%%%%%%%%%%%%%%%%%%%%%%%%%%%%%%%%%%%%%%%%%%%%%%%%%%%%%%%%%%%%%%%%%%%%%%%%%%%%%%%%%%%%%%%%%%%%%
%%%%%%%%%%%%%%%%%%%%%%%%%%%%%%%%%%%%%%%%%%%%%%%%%%%%%%%%%%%%%%%%%%%%%%%%%%%%%%%%%%%%%%%%%%%%%%%%%%%%%%%%%%%%%%%%%%%%%%%%%%%%%%%%%%%%%%%%%%%%%%%%%%%%%%%%%%%%

%%%%%%%%%%%%%%%%%%%%%%%%%%%%%%%%%%%%%%%%%%%%%%%%%%%%%%%%%%%%%%%%%%%%%%%%%%%%%%%%%%%%%%%%%%%%%%%%%%%%%%%%%%%%%%%%%%%%%%%%%%%%%%%%%%%%%%%%%%%%%%%%%%%%%%%%%%%%%%%%%%

\begin{figure}
\begin{center}
\raisebox{0ex}{\includegraphics[width=0.45\linewidth]{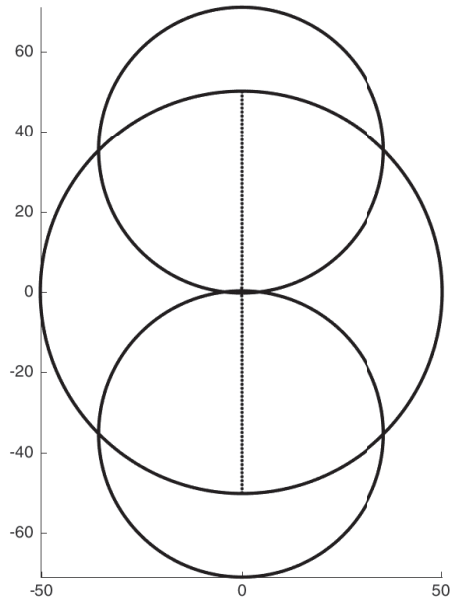}}
\caption{Inclusion regions for Example~3.}
\label{bgb_fig9}                      
\end{center}
\end{figure}

%%%%%%%%%%%%%%%%%%%%%%%%%%%%%%%%%%%%%%%%%%%%%%%%%%%%%%%%%%%%%%%%%%%%%%%%%%%%%%%%%%%%%%%%%%%%%%%%%%%%%%%%%%%%%%%%%%%%%%%%%%%%%%%%%%%%%%%%%%%%%%%%%%%%%%%%%%%%%%%%%%
%
%
%                                                         EXAMPLE 3                                                               
%
%
%%%%%%%%%%%%%%%%%%%%%%%%%%%%%%%%%%%%%%%%%%%%%%%%%%%%%%%%%%%%%%%%%%%%%%%%%%%%%%%%%%%%%%%%%%%%%%%%%%%%%%%%%%%%%%%%%%%%%%%%%%%%%%%%%%%%%%%%%%%%%%%%%%%%%%%%%%%%%%%%%%
\noindent {\bf Example 3.} This example is taken from~\cite{HighamTisseur}. Its quadratic matrix polynomial $Iz^{2} + \Aone z + \Azero$ originates
from a Galerkin method with $n$ basis functions applied to a second order partial differential equation
describing the free vibration of a string, clamped at both ends in a spatially inhomogeneous environment. Here the matrix coefficients are given, for 
$\epsilon,\delta > 0$, by
\bdis
A_{0} = \pi \underset{1 \leq j \leq n}{\operatorname{diag}} \lb j^{2} \rb 
\;\; , \;\; 
\lb A_{1} \rb_{k\ell} = 2 \epsilon \int_{0}^{\pi} \lb x^{2}(\pi - x)^{2} - \delta \rb \sin{(kx)}\sin{(\ell x)} \, dx \; .
\edis
With $n=50$, $\epsilon=0.1$, and $\delta=2.7$ as in~\cite{HighamTisseur},
we proceed as in the previous example, using similar arguments for the choice of the nodes. As for the previous example, we have shown results only for 
the Newton basis as there is very little difference in the size of the disks between the $\mathcal{N}$ and $\mathcal{B}$ bases. 
Figure~\ref{bgb_fig9} shows the inclusion region for the Newton basis.
The eigenvalues are concentrated along the imaginary axis, but they are not purely imaginary.
%
% Figures created with bgb_example3.m 
%
%%%%%%%%%%%%%%%%%%%%%%%%%%%%%%%%%%%%%%%%%%%%%%%%%%%%%%%%%%%%%%%%%%%%%%%%%%%%%%%%%%%%%%%%%%%%%%%%%%%%%%%%%%%%%%%%%%%%%%%%%%%%%%%%%%%%%%%%%%%%%%%%%%%%%%%%%%%%
%%%%%%%%%%%%%%%%%%%%%%%%%%%%%%%%%%%%%%%%%%%%%%%%%%%%%%%%%%%%%%%%%%%%%%%%%%%%%%%%%%%%%%%%%%%%%%%%%%%%%%%%%%%%%%%%%%%%%%%%%%%%%%%%%%%%%%%%%%%%%%%%%%%%%%%%%%%%
%
%                                           F I G U R E  9
%
%%%%%%%%%%%%%%%%%%%%%%%%%%%%%%%%%%%%%%%%%%%%%%%%%%%%%%%%%%%%%%%%%%%%%%%%%%%%%%%%%%%%%%%%%%%%%%%%%%%%%%%%%%%%%%%%%%%%%%%%%%%%%%%%%%%%%%%%%%%%%%%%%%%%%%%%%%%%
%%%%%%%%%%%%%%%%%%%%%%%%%%%%%%%%%%%%%%%%%%%%%%%%%%%%%%%%%%%%%%%%%%%%%%%%%%%%%%%%%%%%%%%%%%%%%%%%%%%%%%%%%%%%%%%%%%%%%%%%%%%%%%%%%%%%%%%%%%%%%%%%%%%%%%%%%%%%
%%%%%%%%%%%%%%%%%%%%%%%%%%%%%%%%%%%%%%%%%%%%%%%%%%%%%%%%%%%%%%%%%%%%%%%%%%%%%%%%%%%%%%%%%%%%%%%%%%%%%%%%%%%%%%%%%%%%%%%%%%%%%%%%%%%%%%%%%%%%%%%%%%%%%%%%%%%%
%%%%%%%%%%%%%%%%%%%%%%%%%%%%%%%%%%%%%%%%%%%%%%%%%%%%%%%%%%%%%%%%%%%%%%%%%%%%%%%%%%%%%%%%%%%%%%%%%%%%%%%%%%%%%%%%%%%%%%%%%%%%%%%%%%%%%%%%%%%%%%%%%%%%%%%%%%%%

\noindent {\bf Summary}
\newline
We have derived inclusion regions for the eigenvalues of matrix polynomials expressed in a general basis and shown the advantages this can provide at the 
hand of several examples from the engineering literature. Not every problem benefits from a change of basis, but there is apparently no shortage of problems
that do. We further remark that the relatively crude estimations we have used to determine the nodes of the bases $\mathcal{N}$ and $\mathcal{B}$
will generally be different for different problems and may be refined, depending on the properties of the coefficient matrices and the choice of matrix
norm. 
Fortunately, the computational cost involved is negligible compared to the computation of the eigenvalues themselves, so that there is no reason not
to try and use a more general basis, especially since the eigenvalues must lie in the intersection of all the inclusion regions obtained for different bases,
further reducing the size of those regions. Finally, we mention that the reverse polynomial can be used to generate additional information on the location of the 
eigenvalues.


\begin{thebibliography}{10}

\bibitem{BHMST}
Betcke, T., Higham, N.J., Mehrmann, V., Schröder, C., and Tisseur, F.
\emph{NLEVP: a collection of nonlinear eigenvalue problems.}
\newblock ACM Trans. Math. Software, 39 (2013), no. 2, Art. 7, 28~pp. 

\bibitem{BiniNoferiniSharify}
Bini, D.A., Noferini, V., and Sharify, M.
\emph{Locating the eigenvalues of matrix polynomials.}     
\newblock SIAM J. Matrix Anal. Appl., 34 (2013), 1708--1727.

\bibitem{Cauchy}
Cauchy, A.L.
\emph{Sur la r\'{e}solution des \'{e}quations num\'{e}riques et sur la th\'{e}orie de l'\'{e}limination.}
\newblock Exercices de Math\'{e}matiques, Quatri\`{e}me Ann\'{e}e, p.65--128. de Bure fr\`{e}res, Paris, 1829.
\newblock Also in: Oeuvres Compl\`{e}tes, S\'{e}rie 2, Tome 9, 86--161.
Gauthiers-Villars et fils, Paris, 1891.

\bibitem{CCVG}
Chaitin-Chatelin, F. and van Gijzen, M.B.
\emph{Analysis of parameterized quadratic eigenvalue problems in computational acoustics with homotopic deviation theory.}
\newblock Numer. Linear Algebra Appl., 13 (2006), 487--512.

%\bibitem{FPS}
%Feriani, A., Perotti, F., and Simoncini, V.
%\emph{Iterative system solvers for the frequency analysis of linear mechanical systems.}
%\newblock Comput. Methods Appl. Mech, Eng., 19 (2000), 1719--1739.

\bibitem{HT_pseudospectra} 
Higham, N.J and Tisseur, F.
\emph{More on pseudospectra for polynomial eigenvalue problems and applications in control theory.}
\newblock Fourth special issue on linear systems and control. Linear Algebra Appl., 351/352 (2002), 435--453.

\bibitem{HighamTisseur}
Higham, N.J. and Tisseur, F.
\emph{Bounds for eigenvalues of matrix polynomials.}
\newblock Linear Algebra Appl., 358 (2003), 5--22. 

\bibitem{HJ}
Horn, R. A. and Johnson, C. R.
\emph{Matrix Analysis.}
\newblock Cambridge University Press, Cambridge, 2013.

\bibitem{Marden}          
Marden, M.
\emph{Geometry of polynomials.}
\newblock Mathematical Surveys, No. 3, American Mathematical Society, Providence, R.I., 1966.

\bibitem{Melman_MatPol}             
Melman, A.
\emph{Generalization and variations of Pellet's theorem for matrix polynomials.}
\newblock Linear Algebra Appl., 439 (2013), 1550–1567. 

\bibitem{RS}
Rahman, Q.I., and Schmeisser, G.
\emph{Analytic Theory of Polynomials}
\newblock London Mathematical Society Monographs. New Series, 26. The Clarendon Press, Oxford University Press, Oxford, 2002.

\bibitem{SP}
Simoncini, V. and Perotti, F.
\emph{On the numerical solution of ($\lambda^2 A+\lambda B+C)x=b$ and application to structural dynamics.}
\newblock SIAM J. Sci. Comput., 23 (2002), 1875–1897.

\bibitem{TisseurMeerbergen}
Tisseur, F. and Meerbergen, K.
\emph{The quadratic eigenvalue problem.}
\newblock SIAM Rev., 43 (2001), 235--286.  

\bibitem{TH_pseudospectra}
Tisseur, F. and Higham, N.J.
\emph{Structured pseudospectra for polynomial eigenvalue problems, with applications.}
\newblock SIAM J. Matrix Anal. Appl., 23 (2001), 187--208 (electronic). 

\end{thebibliography}
\end{document}